\documentclass[a4paper,10pt]{amsart}
\usepackage{latexsym}
\usepackage{graphicx}

\def\oop#1{}

\def\PicV#1#2{{\includegraphics[scale=#1]{#2.eps}}}

\def\JustShowFigureV#1#2{\vskip5pt\begin{center}\PicV{#1}{#2}\end{center}\vskip5pt}

\def\ShowFigureVI#1#2#3{\vskip5pt\noindent \begin{center}\PicV{#1}{#2}\vskip2pt\centerline{Figure
\thefigcounter: #3}\end{center}\newcounter{#2}\setcounter{#2}{\thefigcounter}
\addtocounter{figcounter}{1}\vskip5pt}

\def\ShowFigureVE#1#2#3#4#5#6{\begin{center}
\begin{tabular}{lcccr}
\PicV{#1}{#3} & \,\,\,\,\, & \PicV{#1}{#4}& \,\,\,\,\, &
\PicV{#1}{#5}\\[-2pt]
{\small (a)}&&{\small (b)}&&{\small (c)}\\[-2pt]
\end{tabular}\\[6pt]
Figure \thefigcounter: #6
\newcounter{#2}\setcounter{#2}{\thefigcounter}
\addtocounter{figcounter}{1}
\end{center}}

\input xy
\xyoption{all}

\CompileMatrices

\newtheorem{theorem}{Theorem}

\newtheorem{lemma}[theorem]{Lemma}
\newtheorem{corollary}[theorem]{Corollary}

\DeclareMathOperator{\Aut}{Aut} %
\DeclareMathOperator{\Mon}{Mon} %
\DeclareMathOperator{\Stab}{Stab} %

\DeclareMathOperator{\Alg}{Alg} %
\def\cG{\mathcal{G}}
\def\cH{\mathcal{H}}
\def\cR{\mathcal{R}}

\def\cM{\mathcal{M}}

\def\cD{\mathcal{D}}
\def\cK{\mathcal{K}}
\def\cS{\mathcal{S}}
\def\cP{\mathcal{P}}
\def\cQ{\mathcal{Q}}
\def\cC{\mathcal{C}}
\def\cT{\mathcal{T}}
\def\cI{\mathcal{I}}
\def\cO{\mathcal{O}}

\def\De{\Delta}

\def\mN{\mathbb{N}}
\def\isomorphic{\cong}
\def\st{\mid}

\def\para{\longrightarrow}
\def\mo{\!\!^{-1}}
\def\entao{\Rightarrow}
\def\se{\Leftarrow}
\def\sse{\Leftrightarrow}
\def\x{\times}

\def\vphi{\varphi}
\def\normal{\lhd}
\def\V{\forall}
\def\nn{\cap}
\def\cc{\subset}
\def\m{\backslash}
\def\gpg#1{\langle #1\rangle}
\def\<{\langle }
\def\>{\rangle }

\def\sub#1{_{_{#1}}}
\def\sup#1{^{^{#1}}}

\def\rco{/\!\! _{_r}}

\newcommand{\nclosure}[2]{{#1^{^{#2}}}}

\def\nclosureD#1{\langle #1\rangle^{\Delta}}
\def\nclosureG#1{\langle #1\rangle^{G}}

\newcommand{\core}[2]{#1_{_{#2}}}

\def\coreD#1{{#1}\sub{\Delta}}
\def\coreDho#1{{#1}_{\Delta^{\hat 0}}}
\def\ncD#1{{#1}\sup{\Delta}}

\def\cDH{{H\sub{\Delta}}}
\def\ccover#1{{#1}^{\Delta}}
\def\ccore#1{{#1}_{\Delta}}
\def\Dho{\Delta^{\hat 0}}
\def\Dhi{\Delta^{\hat 1}}
\def\Dhz{\Delta^{\hat 2}}
\def\Do{\Delta^{0}}
\def\Di{\Delta^{1}}
\def\Dz{\Delta^{2}}
\def\Dp{\Delta^{+}}
\def\Dpho{\Delta^{+\hat 0}}
\def\Dphz{\Delta^{+\hat 2}}

\def\hide#1{}
\def\ifif{if and only if }
\def\tx#1{\mbox{#1}}

\def\noi{\noindent}

\def\tfH{{\cH}^{{\vphi}\mo}}

\def\Wal{Wal}
\def\Pin{Pin}
\def\fw{\varphi_{_W}}
\def\fp{\varphi_{_P}}
\def\op{^{\!\circ}}

\def\ig{\Upsilon}

\def\uigD{\Upsilon\sup{\Delta}}
\def\ligD{\Upsilon\sub{\Delta}}

\def\ii{\iota}

\def\uiiD{\iota\sup{\Delta}}
\def\liiD{\iota\sub{\Delta}}

\def\chS{$\check{{\rm S}}$}

\begin{document}
\newcounter{figcounter}
\addtocounter{figcounter}{1}

\title[Bipartite-uniform hypermaps on the sphere]{Bipartite-uniform hypermaps on the sphere}

\author{Antonio Breda d'Azevedo \and Rui Duarte}
\address{Department of Mathematics \\
University of Aveiro \\ 3810-193 Aveiro \\ Portugal}

\thanks{Research partially supported by R\&DU ``Matem\'{a}tica e Aplica\c{c}\~{o}es'' of
Universidade de Aveiro through ``Programa Operacional Ci\^{e}ncia, Tecnologia, Inova\c{c}\~{a}o'' (POCTI) of
the ``Funda\c{c}\~{a}o para a Ci\^{e}ncia e a Tecnologia'' (FCT), cofinanced by the European Community fund
FEDER.}

\date{May 28, 2006}

\begin{abstract}
A hypermap is (hypervertex-) bipartite if its hypervertices can be 2-coloured in such a way that
``neighbouring'' hypervertices have different colours. It is bipartite-uniform if within each of
the sets of hypervertices of the same colour, hyperedges and hyperfaces, elements have common
valencies. The flags of a bipartite hypermap are naturally 2-coloured by assigning the colour of
its adjacent hypervertices. A hypermap is bipartite-regular if the automorphism group acts
transitively on each set of coloured flags. If the automorphism group acts transitively on the set
of all flags, the hypermap is regular. In this paper we classify the bipartite-uniform hypermaps
on the sphere (up to duality). Two constructions of bipartite-uniform hypermaps are given. All
bipartite-uniform spherical hypermaps are shown to be constructed in this way. As a by-product we
show that every bipartite-uniform hypermap $\cH$ on the sphere is bipartite-regular. We also
compute their irregularity group and index, and also their closure cover $\ccover{\cH}$ and
covering core $\ccore{\cH}$.
\end{abstract}

\maketitle

\section{Introduction}

A map generalises to a hypermap when we\oop{waive} lay off the requirement that an edge must join
two vertices at most. A hypermap $\cH$ can be regarded as a bipartite map where the two
monochromatic sets of vertices represent one, the hypervertices and the other the hyperedges of
$\cH$. In this perspective hypermaps are cellular embeddings of hypergraphs on compact connected
surfaces (two-dimensional compact connected manifolds) without boundary $-$ in this paper we deal
only with the boundary-free case.

Usually classifications in map/hypermap theory are carried out by genus, by number of faces, by
embedding of graphs, by automorphism groups or by some fixed properties such as edge-transitivity.
Since Klein and Dyck \cite{Klein,Dyck1} -- where certain 3-valent regular maps of genus 3 were
studied in connection with constructions of automorphic functions on surfaces -- that most
classifications on maps (and hypermaps) involve regularity or orientably-regularity
(direct-regularity). The orientably-regular maps on the Torus (in \cite{CM}), the
orientably-regular embeddings of complete graphs (in \cite{JJroicg}), the orientably-regular maps
with automorphism groups isomorphic to $PSL(2,q)$ (in \cite{Sgrcrs}) and the bicontactual regular
maps (in \cite{Wbrm}), are examples to name but a few. The just-edge-transitive maps of Jones
\cite{Jjetmcg} and the classification by Siran, Tucker and Watkins \cite{STWrfetom} of the
edge-transitive maps on the Torus, on the other hand, include another kind of ``regularity'' other
than regularity or orientably-regularity. According to Graver and Wakins \cite{GWlfpetg}, an edge
transitive map is determined by 14 types of automorphism groups. Among these, 11 correspond to
``restricted regularity'' \cite{Btrrh}. The Jones's ``just-edge-transitive'' maps correspond to
$\Delta^{\hat 0\hat 2}$-regular maps of ``rank 4'', where $\Delta^{\hat 0\hat 2}$ is the normal
closure of $\gpg{R_1, R_0R_2}$ with index 4 in the free product $\De=C_2*C_2*C_2$ generated by the
3 reflections $R_0$, $R_1$ and $R_2$ on the sides of a hyperbolic triangle with zero internal
angles; ``rank 4'' means that it is not $\Theta$-regular for any normal subgroup of $\Delta$ of
index $<4$.  Moreover, the automorphism group of the toroidal edge-transitive maps realise 7 of
the above 14 family-types \cite{STWrfetom}; they all correspond to restrictedly regular maps,
namely of ranks 1 [the regular maps], 2 [the just-orientably-regular (or chiral) maps, the
just-bipartite-regular maps, the just-face-bipartite-regular maps and the
just-Petrie-bipartite-regular maps] and 4 [the just-$\Dpho$-regular maps and the
just-$\Dphz$-regular maps] (see \cite{Btrrh}).

In this paper we classify the ``bipartite-uniform'' hypermaps on the sphere. They all emerged as
being ``bipartite-regular''. A hypermap $\cH$ is {\em bipartite} if its hypervertices can be
2-coloured in such a way that ``neighbouring'' hypervertices have different colours. It is {\em
bipartite-uniform} if the hypervertices of one colour, the hypervertices of the other colour, the
hyperedges and the hyperfaces have common valencies $l_1$, $l_2$, $m$ and $n$ respectively. Flags
in a bipartite hypermap are naturally 2-coloured by assigning the colour of their adjacent
hypervertices. A bipartite hypermap is {\em bipartite-regular} if the automorphism group acts
transitively on each set of coloured flags. If the automorphism group acts transitively on the
whole set of flags the hypermap is {\em regular}. Bipartite-regularity corresponds to
$\Dho$-regularity \cite{Btrrh} where $\Dho$, a normal subgroup of index 2 in $\De$, is the normal
closure of the subgroup generated by $R_1$ and $R_2$.

 We also compute the irregularity group and the
irregularity index of the bipartite-regular hypermaps $\cH$ on the sphere as well as their closure
cover $\ccover{\cH}$ (the smallest regular hypermap that covers $\cH$) and their covering core
$\ccore{\cH}$ (the largest regular hypermap covered by $\cH$). Regular hypermaps on the sphere are
well known and their classification is an academic exercise. Also well known is the following
fact, which comes from the ``universality'' of the sphere: uniform hypermaps on the sphere are
regular. According to \cite{Btrrh} this translates to ``$\Delta$-uniformity in the sphere implies
$\Delta$-regularity''. We may now ask for which normal subgroups $\Theta$ of finite index in
$\Delta$ do we still have ``$\Theta$-uniformity in the sphere implies $\Theta$-regularity'', once
the meaning of $\Theta$-uniformity is understood? As a byproduct of the classification we show in
this paper that bipartite-uniformity (that is, $\Dho$-uniformity) still implies
bipartite-regularity (that is, $\Dho$-regularity). $\Dho$ is just one of the seven normal
subgroups with index 2 in $\Delta$. The others are $\Dhi=\nclosureD{R_0,R_2}$,
$\Dhz=\nclosureD{R_0,R_1}$, $\Do=\nclosureD{R_0,R_1R_2}$, $\Di=\nclosureD{R_1,R_0R_2}$,
$\Dz=\nclosureD{R_2,R_0R_1}$ and $\Dp=\gpg{R_1R_2,R_2R_0}$ (see \cite{BJdcrah} for more details).
As the notation indicates they are grouped into three families, within which they differ by a dual
operation. This duality says that the result is still valid if we replace $\Dho$ by $\Dhi$ or
$\Dhz$. For $\Theta=\Do,\Di,\Dz$, and $\Dp$, $\Theta$-uniformity is the same as uniformity, and
since regularity implies $\Theta$-regularity, on the sphere $\Theta$-uniformity implies
$\Theta$-regularity for any subgroup $\Theta$ of index 2 in $\Delta$. At the end, as a final
comment, we show that on each orientable surface we can find always bipartite-chiral (that is,
irregular bipartite-regular) hypermaps.

\subsection{Hypermaps}

A \emph{hypermap} is combinatorially described by a four-tuple $\cH = (\Omega_{\cH}; h_0, h_1,
h_2)$ where $\Omega_{\cH}$ is a non-empty finite set and $h_0, h_1, h_2$ are fixed-point free
involutory permutations of $\Omega_{\cH}$ generating a permutation group $\langle h_0, h_1, h_2
\rangle$ acting transitively on $\Omega_{\cH}$. The elements of $\Omega_{\cH}$ are called
\emph{flags}, the permutations $h_0$, $h_1$ and $h_2$ are called \emph{canonical generators} and
the group $\Mon (\cH) = \langle h_0, h_1, h_2 \rangle$ is the \emph{monodromy group of $\cH$}. One
says that $\cH$ is a \emph{map} if $(h_0 h_2)^2 = 1$. The \emph{hypervertices} (or
\emph{$0$-faces}) of $\cH$ correspond to $\langle h_1, h_2 \rangle$-orbits on $\Omega_{\cH}$.
Likewise, the \emph{hyperedges} (or \emph{$1$-faces}) and \emph{hyperfaces} (or \emph{$2$-faces})
correspond to $\langle h_0, h_2 \rangle$ and $\langle h_0, h_1 \rangle$-orbits on $\Omega_{\cH}$,
respectively. If a flag $\omega$ belongs to the corresponding orbit determining a $k$-face $f$ we
say that $\omega$ belongs to $f$, or that $f$ contains $\omega$.

We fix $\{ i, j, k \} = \{ 0, 1, 2 \}$. The \emph{valency} of a $k$-face $f = w\langle h_i, h_j
\rangle$, where $\omega \in \Omega_{\cH}$, is the least positive integer $n$ such that $(h_i
h_j)^n \in \Stab (w)$. Since $h_i \neq 1 $ and $h_j\neq 1$, $h_i h_j$ generates a normal subgroup
with index two in $\langle h_i, h_j \rangle$. It follows that $|\langle h_i, h_j \rangle| = 2
|\langle h_i h_j \rangle|$ and so the valency of a $k$-face is equal to half of its cardinality.
$\cH$ is \emph{uniform} if its $k$-faces have the same valency $n_k$, for each $k \in \{ 0, 1, 2
\}$. We say that $\cH$ has \emph{type} $(l; m; n)$ if $l$, $m$ and $n$ are, respectively, the
least common multiples of the valencies of the hypervertices, hyperedges and hyperfaces. The
\emph{characteristic} of a hypermap is the Euler characteristic of its underlying surface, the
imbedding surface of the underlying hypergraph.

A \emph{covering} from a hypermap $\cH = (\Omega_{\cH}; h_0, h_1, h_2)$ to another hypermap $\cG =
(\Omega_{\cG}; g_0, g_1, g_2)$ is a function $\psi: \Omega_{\cH} \to \Omega_{\cG}$ such that $h_i
\psi = \psi g_i$ for all $i \in \{ 0, 1, 2 \}$. The transitive action of $\Mon (\cG)$ on
$\Omega_{\cG}$ implies that $\psi$ is onto. By von Dyck's theorem (\cite{Johnson}, pg 28) the
assignment $h_i \mapsto g_i$ extends to a group epimorphism $\Psi: \Mon (\cH) \to \Mon (\cG)$
called the \emph{canonical epimorphism}. The covering $\psi$ is an \emph{isomorphism} if it is
injective. If there exists a covering $\psi$ from $\cH$ to $\cG$, we say that $\cH$ \emph{covers}
$\cG$ or that $\cG$ \emph{is covered by} $\cH$; if $\psi$ is an isomorphism we say that $\cH$ and
$\cG$ are \emph{isomorphic} and write $\cH \cong \cG$. An \emph{automorphism} of $\cH$ is an
isomorphism $\psi: \Omega_{\cH}\to \Omega_{\cH}$ from $\cH$ to itself; that is, a function $\psi$
that commutes with the canonical generators. The set of automorphisms of $\cH$ is represented by
$\Aut (\cH)$. As a direct consequence of the Euclidean Division Algorithm we have:

\begin{lemma}
Let $\psi: \Omega_{\cH} \to \Omega_{\cG}$ be a covering from $\cH$ to $\cG$ and $\omega \in
\Omega_{\cH}$. Then the valency of the $k$-face of $\cG$ that contains $\omega \psi$ divides the
valency of the $k$-face of $\cH$ that contains $\omega$.
\end{lemma}

Of the two groups $\Mon(\cH)$ and $\Aut(\cH)$ the first acts transitively on $\Omega=\Omega_{\cH}$
(by definition) and the second, due to the commutativity of the automorphisms with the canonical
generators, acts {\em semi-regularly} on $\Omega$; that is, the non-identity elements of
$\Aut(\cH)$ act without fixed points. A transitive semi-regular action is called a {\em regular}
action. These two actions give rise to the following inequalities:

\[|\Mon (\cH)| \geq |\Omega| \geq |\Aut (\cH)|\,.\]

\noindent Moreover, each of the above equalities implies the other. An  equality in the first
inequalities implies that $Mon(\cH)$ acts semi-regularly (hence regularly) on $\Omega$, while an
equality on the the second implies that $Aut(\cH)$ acts transitively (hence regularly) on
$\Omega$. If $\Mon (\cH)$ acts regularly on $\Omega$, or equivalently if $\Aut (\cH)$ acts
regularly on $\Omega$, the hypermap $\cH$ is regular.

\vskip4pt

Each hypermap $\cH$ gives rise to a permutation representation $\rho_{\cH} : \Delta \to \Mon
(\cH)$, $R_i \mapsto h_i$, where $\Delta$ is the free product $C_2*C_2*C_2$ with presentation
$\Delta = \langle R_0, R_1, R_2\mid {R_0}^2 = {R_1}^2 = {R_2}^2 = 1 \rangle$. The group $\Delta$
acts naturally and transitively on $\Omega_{\cH}$ via $\rho_{\cH}$. The stabiliser $H =
\Stab_{\Delta} (\omega)$ of a flag $\omega \in \Omega_{\cH}$ under the action of $\Delta$ is
called the \emph{hypermap subgroup} of $\cH$; this is unique up to conjugation in $\Delta$. The
valency of a $k$-face containing $\omega$ is the least positive integer $n$ such that $(R_i R_j)^n
\in H$; more generally, the valency of a $k$-face containing the flag $\sigma = \omega\cdot g =
\omega (g) \rho_{\cH} \in \Omega_{\cH}$, where $g \in \Delta$, is the least positive integer $n$
such that $(R_i R_j)^n \in \Stab_{\Delta} (\sigma) = \Stab_{\Delta} (\omega\cdot g) =
{\Stab_{\Delta} (\omega)}^g = H^g$. \vskip4pt Denote by $\Alg (\cH) = (\Delta /_{\small r} H; a_0,
a_1, a_2)$ where $a_i: \Delta /_{\small r} H \to \Delta /_{\small r} H$, $Hg \mapsto Hg \cDH R_i =
HgR_i$. It is easy to see that $\Alg (\cH) \cong \cH$. We say that $\Alg (\cH)$ is the
\emph{algebraic presentation} of $\cH$. Moreover, it is well known that:

\begin{enumerate}
\item A hypermap $\cH$ is regular if and only if its hypermap subgroup $H$ is normal in $\Delta$.
\item A regular hypermap is necessarily uniform.
\end{enumerate}

\noindent Since $\Alg (\cH)$ and $\cH$ are isomorphic, we will not differentiate one from the
other. \vskip5pt

Following \cite{Btrrh} if $H < \Theta$ for a given $\Theta \normal \Delta$, we say that $\cH$ is
\emph{$\Theta$-conservative}. A $\Delta^+$-conservative hypermap is better known as an orientable
hypermap. An automorphism of an orientable hypermap either preserves $\Dp$-orbits or not. Those
that preserve $\Dp$-orbits are called {\em orientation-preserving automorphisms}. The set of
orientation-preserving automorphisms is a subgroup of $\Aut(\cH)$ and is denoted by $\Aut^+(\cH)$.
If $\cH$ is $\Delta^{\hat 0}$-conservative (resp. $\Delta^{\hat 1}$-conservative, resp.
$\Delta^{\hat 2}$-conservative) we say that $\cH$ is \emph{bipartite}, \emph{vertex-bipartite} or
\emph{$0$-bipartite} (resp. \emph{edge-bipartite} or \emph{$1$-bipartite}, resp.
\emph{face-bipartite} or \emph{$2$-bipartite}).

\begin{lemma} \label{valenciaspares}
If $\cH$ is bipartite and $\omega \in \Omega_{\cH}$, then the valencies of the hyperedge and the
hyperface that contain $\omega$ must be even.
\end{lemma}
\begin{proof}
If $m$ and $n$ are the valencies of the hyperedge and the hyperface that contain $\omega=Hd$,
$d\in\De$, then $(R_2 R_0)^{m}, (R_0 R_1)^{n} \in H^d \subseteq \Dho$. Therefore $m$ and $n$ must
be even.
\end{proof}

If $H \lhd \Delta^+$, we say that $\cH$ is \emph{orientably-regular}. If $H \lhd \Delta^{\hat 0}$
(resp. $H \lhd \Delta^{\hat 1}$ and $H \lhd \Delta^{\hat 2}$), we say that $\cH$ is
\emph{vertex-bipartite-regular} (resp. \emph{edge-bipartite-regular} and
\emph{face-bipartite-regular}). If $\cH$ is vertex-bipartite-regular (resp.
edge-bipartite-regular, resp. face-bipartite-regular) but not regular, we say that $\cH$ is
\emph{vertex-bipartite-chiral} (resp. \emph{edge-bipartite-chiral}, resp.
\emph{face-bipartite-chiral}). We will use \emph{bipartite-regular} and \emph{bipartite-chiral} in
place of vertex-bipartite-regular and vertex-bipartite-chiral for short.

A bipartite-uniform hypermap is a bipartite hypermap such that all the hypervertices in the same
$\Delta^{\hat 0}$-orbit have the same valency, as do all the hyperedges and all the hyperfaces.
The \emph{bipartite-type} of a bipartite-uniform hypermap $\cH$ is a four-tuple $(l_1, l_2; m; n)$
(or $(l_2, l_1; m; n)$) where $l_1$ and $l_2$ ($l_1 \leq l_2$) are the valencies (not necessarily
distinct) of the hypervertices of $\cH$, $m$ is the valency of the hyperedges of $\cH$ and $n$ is
the valency of the hyperfaces of $\cH$. We note that if $\cH$ is a bipartite-uniform hypermap of
bipartite-type $(l_1, l_2; m; n)$, then $m$ and $n$ must be even by Lemma \ref{valenciaspares}.

\subsection{Euler formula for uniform hypermaps}

Using the well known Euler formula for maps one easily gets the following well known result:

\begin{lemma}[Euler formula for hypermaps] \label{effh}
Let $\cH$ be a hypermap with $V$ hypervertices, $E$ hyperedges and $F$ hyperfaces. If $\cH$ has
underlying surface $\cS$ with Euler characteristic $\chi$, then  $\chi = V + E + F -
\frac{|\Omega_{\cH}|}{2}$.
\end{lemma}

If $\cH$ is uniform of type $(l, m, n)$, then $V = \frac{|\Omega_{\cH}|}{2l}$, $E =
\frac{|\Omega_{\cH}|}{2m}$ and $F = \frac{|\Omega_{\cH}|}{2n}$. Replacing the values of $V$, $E$
and $F$ in the last formula, we get:

\begin{corollary}[Euler formula for uniform hypermaps] \label{effuh}
\[\chi = \frac{|\Omega_{\cH}|}{2} \left( \frac{1}{l} + \frac{1}{m} + \frac{1}{n} - 1 \right).\]
\end{corollary}

\def\niaut{\psi}

\subsection{Duality}

A non-inner automorphism $\niaut$ of $\Delta$ (that is, an automorphism not arising from a
conjugation) gives rise to an operation on hypermaps by transforming a hypermap $\cH=(\Delta\rco
H,\cDH R_0,\cDH R_1,\cDH R_2)$, with hypermap-subgroup $H$, into its operation-dual
\[\begin{array}{lcl}
  D_{\niaut}(\cH) & = & (\Delta\rco H\niaut;\coreD{(H\niaut)}R_0,\coreD{(H\niaut)}R_1,\coreD{(H\niaut)}R_2) \\[4pt]
                  & = & (\Delta\rco H\niaut;\cDH\niaut R_0,\cDH\niaut R_1,\cDH\niaut R_2) \\[4pt]
\end{array}\]
with hypermap-subgroup $H\niaut$ (see \cite{LyJohoa,Jomh,JTomooa} for more details). Note that if
$\niaut$ is inner, then $D_{\niaut}(\cH)$ is isomorphic to $\cH$. In particular, each permutation
$\sigma \in S_{\{ 0, 1, 2 \}}\m\{id\}$ induces a non-inner automorphism $\sigma\op:\Delta\para
\Delta$ by assigning $R_i\mapsto R_{i\sigma}$, for $i=0,1,2$.  This automorphism induces an
operation $\cD_{\sigma}$ on hypermaps by assigning the hypermap-subgroup $H$ of $\cH$ to a
hypermap-subgroup $H\sigma\op $. Such an operator transforms each hypermap
$\cH=(\Omega_{\cH};h_0,h_1,h_2)$ into its $\sigma$-dual $D_\sigma(\cH)\isomorphic
(\Omega_{\cH};h_{0\sigma\mo},h_{1\sigma\mo},h_{2\sigma\mo})$. We note that the $k$-faces of $\cH$
are the $k \sigma$-faces of $\cD_{\sigma} (\cH)$.

\begin{lemma}[Properties of $\cD_{\sigma}$]
Let $\cH$, $\cG$ be two hypermaps and $\sigma, \tau \in S_{\{ 0, 1, 2 \}}$. Then (1)\, $\cD_{1}
(\cH) = \cH$, where $1 = id \in S_{\{ 0, 1, 2 \}}$;\,\,\, (2)\, $\cD_{\tau} (\cD_{\sigma} (\cH)) =
\cD_{\sigma\tau} (\cH)$;\,\,\, (3)\, If $\cH$ covers $\cG$, then $\cD_{\sigma} (\cH)$ covers
$\cD_{\sigma} (\cG)$;\,\,\, (4)\, If $\cH \cong \cG$, then $\cD_{\sigma} (\cH) \cong \cD_{\sigma}
(\cG)$;\,\,\, (5)\, If $\cH$ is uniform, then $\cD_{\sigma} (\cH)$ is uniform;\,\,\, (6)\, If
$\cH$ is $k$-bipartite-uniform, then $\cD_{\sigma} (\cH)$ is $k \sigma$-bipartite-uniform;\,\,\,
(7)\, If $\cH$ is regular, then $\cD_{\sigma} (\cH)$ is regular;\,\,\, (8)\, If $\cH$ is
$k$-bipartite-regular, then $\cD_{\sigma} (\cH)$ is $k \sigma$-bipartite-regular;\,\, (9) Both
$\cH$ and $\cD_{\sigma} (\cH)$ have same underlying surface.
\end{lemma}

\subsection{Spherical uniform hypermaps}

A hypermap $\cH$ is \emph{spherical} if its underlying surface is a sphere (i.e if its Euler
characteristic is 2). By taking $l \leq m \leq n$ and  $\chi =2$ in the Euler formula one easily
sees that $l < 3$. A simple analysis to the above inequality drives us to the following table of
possible types (up to duality): \vskip5pt
\begin{table}[h]
\[\begin{array}{lll|ccc|ccc|c}
l & m & n & V & E & F & |\Omega_{\cH}| & \Mon(\cH) & \cH & \Aut^+(\cH)\\
\hline
1 & k & k & k & 1 & 1 & 2k & D_k & \cD_{(02)} (\cD_k) & C_k\\
2 & 2 & k & k & k & 2 & 4k & D_k\x C_2 & \cP_k & C_k\\
2 & 3 & 3 & 6 & 4 & 4 & 24 & S_4 & \cD_{(01)} (\cT) & A_4\\
2 & 3 & 4 & 12 & 8 & 6 & 48 & S_4\x C_2 & \cD_{(01)} (\cC) & S_4\\
2 & 3 & 5 & 30 & 20 & 12 & 120 & A_5\x C_2 & \cD_{(01)} (\cD) & A_5\\
\hline
\end{array}\]
\caption{\label{sphereunif} Possible values (up to duality) for type $(l; m;n)$.}
\end{table}

\begin{lemma}
All uniform hypermaps on the sphere are regular.
\end{lemma}

This result arises because each type $(l;m;n)$ in Table \ref{sphereunif} determines a cocompact
subgroup $H=\nclosureD{(R_1R_2)^l, (R_2R_0)^m,(R_0R_1)^n}$ with index $|\Omega_{\cH}|$ in the free
product $\Delta=C_2*C_2*C_2$ generated by $R_0$, $R_1$ and $R_2$.

Let $\cT$, $\cC$, $\cO$, $\cD$ and $\cI$ denote the 2-skeletons of the tetrahedron, the cube, the
octahedron, the dodecahedron and the icosahedron. These are, up to isomorphism, the unique uniform
hypermaps of type $(3;2;3)$, $(3;2;4)$, $(4;2;3)$, $(3;2;5)$ and $(5;2;3)$ respectively, on the
sphere; note that $\cO \cong \cD_{(02)} (\cC)$ and $\cI \cong \cD_{(02)} (\cD)$. Together with the
infinite families of hypermaps $\cD_n$ with monodromy group $D_n$ and $\cP_n$ with monodromy group
$D_n \times C_2$ ($n \in \mN$), of types $(n;n;1)$ and $(2;2;n)$, respectively, they complete, up
to duality and isomorphism, the uniform spherical hypermaps.

\JustShowFigureV{1}{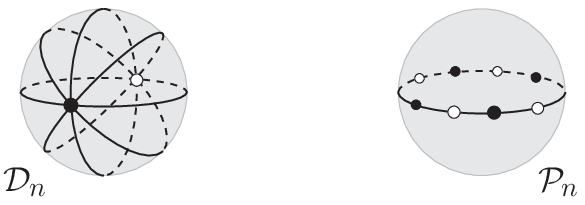}

\noi The last column of Table \ref{sphereunif} displays the uniform spherical hypermaps (which are
regular by last lemma) of type $(l;m;n)$ with $l \leq m \leq n$.

\begin{lemma}\label{valencyone}
If $\cH$ is a hypermap such that all hyperfaces have valency 1, then $\cH$ is the ``dihedral''
hypermap $\cD_n$, a regular hypermap on the sphere with $n$ hyperfaces.
\end{lemma}

\begin{proof}
Let $H$ be a hypermap-subgroup of $\cH$. All hyperfaces having valency 1 implies that $R_0R_1\in
H^d$ for all $d\in\Delta$ (i.e., $R_0R_1$ stabilises all the flags). Then
$H\gpg{R_1,R_2}=H\gpg{R_0,R_2}=H\gpg{R_0,R_1,R_2}=\Delta/_{\small r} H=\Omega$; that is, $\cH$ has
only one hypervertex and one hyperedge. Hence $\cH \cong \cD_n$, where $n$ is the valency of the
hyperedge and the hyperface of $\cH$.
\end{proof}

\section{Constructing bipartite hypermaps}

By the Reidemeister-Schreier rewriting process \cite{Johnson} it can be shown that
\[\Dho \cong C_2 * C_2 * C_2 * C_2 =
\langle R_1 \rangle * \langle R_2 \rangle * \langle {R_1}^{R_0} \rangle * \langle {R_2}^{R_0}
\rangle\,.\] As a consequence we have an epimorphism $\vphi:\Dho\para \Delta$.

Any such epimorphism $\vphi$ induces a transformation (not an operation) of hypermaps, by
transforming each hypermap $\cH=(\Omega\sub{\cH};h_0,h_1,h_2)$ with hypermap subgroup $H$ into a
hypermap $\tfH=(\Omega;t_0,t_1,t_2)$ with hypermap subgroup $H{\vphi}\mo $.

\[\tfH\left\{
\begin{array}{l}
\hspace{14pt}\xymatrix@R=1pc { \Delta \ar@{=}^{2}[d]&\\&} \\[-10pt]
\left.
\begin{array}{l}
\xymatrix@R=1pc {
\Dho \ar@{ -}[dd] \ar@{->}^{\vphi}[r] & \Delta \ar@{-}[dd]\\
&\\
H\vphi\mo\ar@{->}[r] & H }
\end{array}\right\}{\cH}
\end{array}
\right. \] \vskip4pt \noi Algebraically, $\tfH=(\Delta \rco H{\vphi}\mo ; s_0, s_1, s_2)$ with
$s_i=({H{\vphi}\mo})\sub{\Delta}R_i$ acting on $\Omega=\Delta \rco H{\vphi}\mo $ by right
multiplication. Here
 $(H{\vphi}\mo)\sub{\Delta}$ denotes the core of $H{\vphi}\mo$ in $\Delta$.
In the following lemma we list three elementary, but useful, properties of this transformation
$\vphi$.

\begin{lemma}\label{generalprop}
Let $g\in\Delta$, $W=\coreD{(H\vphi\mo)}w\in \De/\coreD{(H\vphi\mo)}=\Mon(\tfH)$ and $H{\vphi}\mo
g\in \Omega$ be a flag of $\tfH$. Then,
\begin{description}
\item[(1)] If $g\in\Dho$, then $(H\vphi\mo)^g=H^{g\vphi}{\vphi}\mo$. If $g\not\in\Dho$, then
$(H\vphi\mo)^g=\left(H^{(gR_0)\vphi}{\vphi}\mo\right)\sup{R_0}$.
\item[(2)] $\coreDho{(H\vphi\mo)}=\coreD{H}\vphi\mo$ and \, $\coreD{(H\vphi\mo)}=\coreD{H}\vphi\mo\nn
\left(\coreD{H}\vphi\mo\right)\sup{R_0}$.
\item[(3)] $W\in \Stab(H{\vphi}\mo g)\,\sse\,w\in (H\vphi\mo)^g \,\sse \, \left\{\begin{array}{ll}
w\vphi\in H^{g\vphi}\,, & \tx{if $g\in\Dho$} \\[6pt]
w^{R_0}\vphi\in H^{(gR_0)\vphi}\,, & \tx{if $g\not\in \Dho$\,.}
\end{array}\right.$
\noi Moreover, $W\in \Stab(H{\vphi}\mo g)$ implies that $w\in\Dho$.
\end{description}
\end{lemma}

\begin{proof}
(1) If $g\in\Dho$, then $x\in H^{g\vphi}{\vphi}\mo \sse x\vphi\in H^{g\vphi} \sse
(x\vphi)^{(g\vphi)\mo} = (x\vphi)^{g\mo\vphi} = x^{g\mo}\vphi \in H \sse x\in (H\vphi\mo)^g$. If
$g\not\in\Dho$, then $gR_0\in\Dho$ and so
$(H\vphi\mo)^g=\left((H\vphi\mo)^{(gR_0)}\right)\sup{R_0}=\left(H^{(gR_0)\vphi}{\vphi}\mo\right)\sup{R_0}$.
\vskip5pt \noindent (2) Since $\vphi$ is onto, the above item translates into these two results.
\vskip5pt \noindent (3) $W\in \Stab(H{\vphi}\mo g)={\Stab(H \varphi^{-1})}^g\sse w\in
(H\vphi\mo)^g$. Since $H\vphi\mo\normal\Dho$, this implies that $w\in\Dho$.

If $g \in \Delta^{\hat 0}$, then $w \in (H\vphi\mo)^g \stackrel{(1)}{=} H^{g \vphi} \vphi^{-1}\sse
w \vphi \in H^{g \vphi}$.

If $g\not\in\Dho$, then $gR_0\in\Dho$ and so, by above, $w \in (H \varphi^{-1})^g\sse w^{R_0} \in
(H \vphi^{-1})^{gR_0} \sse (w^{R_0})\vphi\in H^{(gR_0)\vphi}$.
\end{proof}

\noi Remark: For simplicity we will not distinguish $W$ from $w$, and so we will see $W$ as a word
on $R_0$, $R_1$ and $R_2$ in $\Delta$ instead of a coset word $\coreD{(H\vphi\mo)}w$. \vskip4pt
\begin{theorem}\label{Dhomonodromy}
If \,\,$\cH\,\isomorphic\, \cG^{{\vphi}\mo}$ for some hypermap $\cG$, then
$\Dho$-$\Mon(\cH)\isomorphic \Mon(\cG)$.
\end{theorem}
\begin{proof}
By Lemma \ref{generalprop}(2) we deduce that
$\Dho$-$\Mon(\cH)=\Dho/\coreDho{H}=\Dho/\coreDho{(G\vphi\mo)}= \Dho/\coreD{G}\vphi\mo\isomorphic
\Delta/\coreD{G}=\Mon(\cG)$.
\end{proof}
\vskip4pt Among many possible canonical epimorphisms $\vphi:\Dho \to \Delta$, there are two that
induce transformations preserving the underlying surface, namely $\fw$ and $\fp$ defined by
\[R_1 \fw = R_1, \quad R_2 \fw = R_2, \quad {R_1}^{R_0} \fw = R_0, \quad {R_2}^{R_0} \fw = R_2,\]
\[R_1 \fp = R_1, \quad R_2 \fp = R_2, \quad {R_1}^{R_0} \fp = R_0, \quad {R_2}^{R_0} \fp = R_0.\]

\vskip4pt \noi Denote by $\Wal(\cH)$ the hypermap $\cH^{{\fw}\mo}$ and by $\Pin (\cH)$ the
hypermap $\cH^{{\fp}\mo}$.  $\Wal(\cH)$ is a map; in fact, since $(R_0R_2)^2={R_2}^{R_0}R_2$ and
$((R_0R_2)^2)^{R_0}= R_2{R_2}^{R_0}$ we have $(R_0R_2)^2\fw = ((R_0R_2)^2)^{R_0}\fw = 1$, and
hence, by Lemma \ref{generalprop}(3), for all $g \in \Delta$, $(R_0 R_2)^2\in \Stab(H {\fw}\mo
g)$.

\vskip4pt Both hypermaps $\Wal (\cH)$ and $\Pin(\cH)$ have the same underlying surface as $\cH$
but while $\Wal (\cH)$ is a map (bipartite map since $H {\fw}\mo\subseteq \Dho$), the well known
Walsh bipartite map of $\cH$ \cite{Whvbm,BJdcrah}, $\Pin (\cH)$ is not necessarily a map.

\ShowFigureVI{.7}{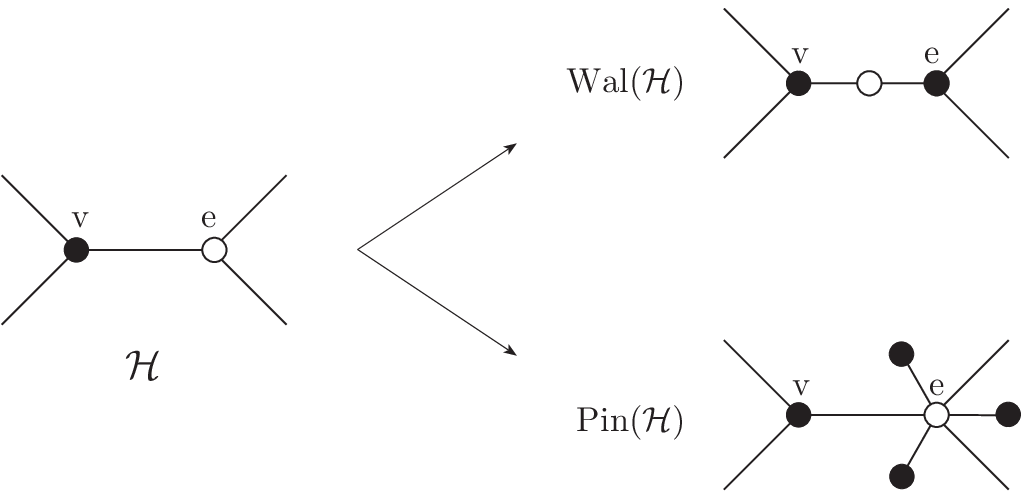}{Topological construction of $\Wal(\cH)$ and $\Pin(\cH)$.}

\begin{theorem}[Properties of $\fw$]\label{pfwalsh} Let $\cH$ be a hypermap. Then:
\begin{enumerate}
\item $\cH$ is uniform of type $(l; m;n)$ \ifif $\Wal(\cH)$ is bipartite-uniform of bipartite-type $(l, m; 2; 2n)$ if $l \leq
m$ or $(m, l; 2; 2n)$ if $l \geq m$;
\item $\cH$ is regular \ifif $\Wal (\cH)$ is bipartite-regular.
\end{enumerate}
\end{theorem}

\begin{proof}
Let $H$ be a hypermap subgroup of $\cH$. Then $H {\fw}\mo $ is a hypermap subgroup of $\Wal
(\cH)$.

 \vskip5pt \noindent (\ref{pfwalsh}.1) ($\entao$) Let us
suppose that $\cH$ is uniform of type $(l; m; n)$. Note first that
\begin{eqnarray}
\label{eqa} R_1 R_2 & = & (R_1 R_2) \fw\,, \\
\label{eqb} R_0 R_2 & = & ({R_1}^{R_0} {R_2}^{R_0}) \fw = (R_1 R_2)^{R_0} \fw\,, \\
\label{eqc} R_0 R_1 & = & ({R_1}^{R_0} R_1) \fw = (R_0 R_1)^2\fw\,.
\end{eqnarray}
Let $W$ denote a word in $R_0$, $R_1$, $R_2$ and $\omega g \in \Omega_{\Wal (\cH)}$ be any flag
($g \in \Delta$). We already know that the valency of the hyperedge containing $\omega g$ is 2
($\Wal(\cH)$ is a map) and that the valency of the hyperface contains $\omega g$ is even. Let $l'$
and $n'$ be the valencies of the hypervertex and the hyperface containing $\omega g$,
respectively. \vskip4pt \noindent (1)\, $g \in \Dho$. From \eqref{eqa} and Lemma
\ref{generalprop}(1) we have $(R_1 R_2)^k \in H^{g \fw}$ if and only if $(R_1 R_2)^{k} \in H^{g
\fw} {\fw}\mo = (H {\fw}\mo)^g$; that is, according to Lemma \ref{generalprop}(3),
\begin{equation}\label{pfweqa}
(R_1R_2)^k\in \Stab(H(g\fw)) \sse (R_1R_2)^k\in \Stab((H\fw\mo) g)\,.
\end{equation}
Analogously, from \eqref{eqc} we get $(R_0 R_1)^k \in H^{g \fw}$ if and only if $(R_0 R_1)^{2k}
\in H^{g \fw} {\fw}\mo = (H {\fw}\mo)^g$ that is, according to Lemma \ref{generalprop}(3),
\begin{equation}\label{pfweqb}
(R_0R_1)^k\in \Stab(H(g\fw)) \sse (R_0R_1)^{2k}\in \Stab((H\fw\mo) g)\,.
\end{equation}
Now the uniformity of $\cH$ implies $l' = l$ and $n' = 2n$.

\vskip4pt \noindent (2)\, $g \notin \Dho$. Since $gR_0\in\Dho$ we get from \eqref{eqb},
\[\begin{array}{rcl}
(R_0 R_2)^k \in H^{(gR_0)\fw} & \sse & ((R_1 R_2)^{R_0})^{k} \in H^{(gR_0)\fw} {\fw}\mo = (H {\fw}\mo)^{gR_0} \\[4pt]
& \sse & (R_1 R_2)^k \in (H {\fw}\mo)^g\,;
\end{array}\]
and from \eqref{eqc},
\[\begin{array}{rcl}
(R_0R_1)^k \in H^{(gR_0)\fw} & \sse & (R_0 R_1)^{2k} \in H^{gR_0\fw}{\fw}\mo = (H{\fw}\mo)^{gR_0} \\[4pt]
& \sse & (R_1R_0)^{2k}\in (H{\fw}\mo)^{g}\,.
\end{array}\]
This implies that

\begin{equation}
\label{pfweqc} (R_0 R_2)^k \in \Stab(H(gR_0)\fw) \sse (R_1 R_2)^k \in \Stab(H\fw\mo g),
\end{equation}
\begin{equation}
\label{pfweqcc} (R_0 R_1)^k \in \Stab(H(gR_0)\fw) \sse (R_1 R_0)^{2k} \in \Stab(H\fw\mo g).
\end{equation}

\noindent Likewise, the uniformity of $\cH$ now implies that $l'= m$ and $n' = 2n$.

\vskip4pt \noindent Gathering (1) and (2) and assuming, without loss of generality, that $l \leq
m$, then $\Wal(\cH)$ is bipartite-uniform of bipartite-type $(l,m;2;2n)$.

\vskip5pt

\noi ($\se$) Let us assume that $\Wal(\cH)$ is bipartite-uniform of bipartite-type $(l,m;2;2n)$.
Being bipartite, $\Wal(\cH)$ has two orbits of vertices: the ``black'' vertices, all with valency
$l$ (say), and the ``white'' vertices, all with valency $m$. Without loss of generality, all the
flags $H\fw\mo g$, $g\in\Dho$, are adjacent to ``black'' vertices while all the flags $H\fw\mo
gR_0$, $g\in\Dho$, are adjacent to ``white'' vertices. As seen before, the equivalence \eqref{eqa}
for $g\in\Dho$ gives rise to the equivalence \eqref{pfweqa}, which expresses the fact that all the
hypervertices of $\cH$ have the same valency $l$; the equivalence \eqref{eqb} for $g\not\in \Dho$
gives rise to the equivalence \eqref{pfweqc}, which says that all the hyperedges of $\cH$ have the
same valency $m$; finally, the equivalence \eqref{eqc} gives rise to the equivalence
\eqref{pfweqb} if $g\in\Dho$ or the equivalence \eqref{pfweqcc} if $g\not\in\Dho$, and they
express the fact that all the hyperfaces of $\cH$ have the same valency $n$. Hence $\cH$ is
uniform of type $(l;m;n)$ (or $(m;l;n)$ since the positional order of $l$ and $m$ in the
bipartite-type of $\Wal(\cH)$ is ordered by increasing value).

\vskip5pt

\noi (\ref{pfwalsh}.2) $\cH \ \text{is regular} \sse H \lhd \Delta \sse H {\fw}\mo\lhd \Dho \sse
\Wal (\cH) \ \text{is bipartite-regular}$ since $\fw$ is an epimorphism.
\end{proof}

\begin{theorem}\label{HWal}
$\cH$ is a bipartite map \ifif $\cH \cong \Wal (\cG)$ for some hypermap $\cG$.
\end{theorem}

\begin{proof}
Only the necessary condition needs to be proved. If $\cH$ is a bipartite map, then $H \subseteq
\Dho$. Since $\cH$ is a map, $((R_0 R_2)^2)^g \in H$ for all $g \in \Delta$; therefore $\ker \fw =
\langle (R_0 R_2)^2 \rangle^{\Dho} \subseteq H$. This implies that $H \fw {\fw}\mo = H \ker \fw =
H$ and hence $\cH \cong \Wal (\cG)$ where $\cG$ is a hypermap with hypermap subgroup $G = H \fw$.
\end{proof}

\begin{theorem}[Properties of $\fp$] Let $\cH$ be a hypermap. Then,
\begin{enumerate}
\item $\Pin (\cH)$ is a bipartite hypermap such that all hypervertices in one $\Dho$-orbit have valency $1$;
\item $\cH$ is uniform of type $(l; m; n)$ \ifif $\Pin (\cH)$ is bipartite-uniform of bipartite-type
$(1, l; 2m; 2n)$;
\item $\cH$ is regular \ifif $\Pin (\cH)$ is bipartite-regular.
\end{enumerate}
\end{theorem}

\begin{proof}
Let $H$ be a hypermap subgroup of $\cH$. Then $H {\fp}\mo $ is a hypermap subgroup of $\Pin
(\cH)$.

\vskip4pt

\noi (1) $\Pin (\cH)$ is bipartite since $H \fp^{-1} \subseteq \Delta^{\hat 0}$. We have $(R_1
R_2)^{R_0} \fp = ({R_1}^{R_0} {R_2}^{R_0}) \fp = 1$; therefore, by Lemma 2 (2), $R_1 R_2 \in
\Stab(H \fp^{-1} g)$ for all $g \not\in \Delta^{\hat 0}$, i.e, all hypervertices in the same
$\Delta^{\hat 0}$-orbit of the hypervertex containing the flag $H \fp^{-1} R_0$ have valency 1.

\vskip4pt

\noi (2) Let us suppose that $\cH$ is uniform of type $(l; m; n)$. We proceed similarly as for
$\fw$, keeping in mind that all hypervertices of $\Pin(\cH)$ adjacent to flags $H\fp\mo g$, for
$g\not\in\Dho$, have valency 1. Starting from the equalities, \noi \[\begin{array}{l}
R_1 R_2 = (R_1 R_2) \fp\,, \\[2pt]
R_0 R_2 = ({R_2}^{R_0} R_2) \fp = (R_0 R_2)^2\fp\,,\\[2pt]
R_0 R_1 = ({R_1}^{R_0} R_1) \fp = (R_0 R_1)^2 \fp\,.
\end{array}\]
one gets the following equivalences,
\[\begin{array}{l}
(R_1R_2)^k\in \Stab(Hg\fp)\sse (R_1R_2)^k\in \Stab(H\fp\mo g),\, \V\, g\in\Dho\,,\\[2pt]
(R_0R_2)^k\in \Stab(Hg\fp)\sse (R_0 R_2)^{2k}\in \Stab(H\fp\mo g),\, \V\, g\in\Dho\,,\\[2pt]
(R_0R_2)^k\in \Stab(H(gR_0)\fp)\sse (R_2 R_0)^{2k}\in \Stab(H\fp\mo g),\, \V\, g\not\in\Dho\,,\\[2pt]
(R_0R_1)^k\in \Stab(Hg\fp)\sse (R_0 R_1)^{2k}\in \Stab(H\fp\mo g),\, \V\, g\in\Dho\,,\\[2pt]
(R_0R_1)^k\in \Stab(H(gR_0)\fp)\sse (R_1 R_0)^{2k}\in \Stab(H\fp\mo g),\, \V\, g\not\in\Dho\,.
\end{array}\]
This clearly shows that $\Pin(\cH)$ is bipartite-uniform of bipartite-type $(1,l;2m;2n)$.
Reciprocally, if $\Pin(\cH)$ is bipartite-uniform of bipartite-type $(1,l;2m;2n)$ then, reversing
the above argument in a similar way as we did for $\Wal(\cH)$ in the proof of Theorem
\ref{pfwalsh}, we easily conclude that $\cH$ is uniform of type $(l;m;n)$.

\vskip4pt

\hide{\noi (3) If $\cH$ is regular, $H \lhd \Delta$ and therefore $H {\fp}\mo\lhd \Dho$, i.e,
$\Pin (\cH)$ is bipartite-regular. Reciprocally, if $\Pin (\cH)$ is bipartite-regular then
$H\fp\mo\normal\Dho$ which implies that $H=H\fp\mo\fp\normal \Dho\fp=\Delta$, that is, $\cH$ is
regular.}

\noi (3) Since $\fp$ is an epimorphism,\,  $\cH \ \text{is regular} \sse H \lhd \Delta \sse H
{\fp}\mo\lhd \Dho \sse \Pin (\cH) \ \text{is bipartite-regular}$.
\end{proof}

\begin{theorem}\label{HPin}
If $\cH$ is a bipartite hypermap such that all hypervertices in one $\Dho$-orbit have valency $1$,
then $\cH \cong \Pin (\cG)$ for some hypermap $\cG$.
\end{theorem}

\begin{proof} As in Theorem 13, only the necessary condition needs to be proved.
Let $H$ be a hypermap subgroup of $\cH$. By taking $H^{R_0}$ instead of $H$ if necessary, we may
assume, without loss of generality, that all hypervertices in the $\Delta^{\hat 0}$-orbit of the
hypervertex that contains the flag $H R_0$ have valency 1, i.e, $R_1 R_2 \in H^{R_0 g}$ for all $g
\in \Delta^{\hat 0}$. Then $((R_1 R_2)^{R_0})^h \in H$ for all $h \in \Delta^{\hat 0}$; therefore
$\ker \fp = \langle (R_1 R_2)^{R_0} \rangle^{\Delta^{\hat 0}} \subseteq H$. This implies that $H
\fp \fp^{\mo} = H \ker \fp = H$ and hence $\cH \cong Pin (\cG)$, where $\cG$ is the hypermap with
hypermap subgroup $G=H\fp$.
\end{proof}

\begin{theorem}
$\Wal (\cD_{(0\,1)} (\cH)) \cong \Wal (\cH)$.
\end{theorem}

\begin{proof}
If $H$ is a hypermap subgroup of $\cH$, then $H {\fw}\mo$ and $H (0\,1)\op {\fw}\mo$ are hypermap
subgroups of $\Wal (\cH)$ and $\Wal (\cD_{(0\,1)} (\cH))$, respectively. Since $g \fw
\sigma=g\iota^{R_0}\fw$ for all $g \in \Dho$, where $\sigma=(0\,1)\op$ and $\iota^{R_0}$ is the
automorphism given by conjugation by $R_0$, we have
\begin{equation}\label{hsigmafwmo}
H\sigma\fw\mo=H\fw\mo\iota^{R_0},
\end{equation}
that is, the hypermap subgroup $H(0\,1)\op\fw\mo$ of $\Wal (\cD_{(0\,1)} (\cH))$ is just a
conjugate under $R_0$ of the hypermap subgroup of $\Wal (\cH)$\oop{,
\[H(0\,1)\op\fw\mo=(H\fw\mo)^{R_0}\,,\]} and so they are isomorphic.
\end{proof}

\begin{theorem}
$\Pin (\cD_{(1\,2)} (\cH)) = \cD_{(1\,2)} (\Pin (\cH))$.
\end{theorem}

\begin{proof}
Let $H$ be a hypermap subgroup of $\cH$ and $\sigma=(1\,2)\op$. Then $H \sigma {\fp}\mo$ and $H
{\fp}\mo \sigma$ are hypermap subgroups of $\Pin (\cD_{(1\,2)} (\cH))$ and $\cD_{(1\,2)} (\Pin
(\cH))$, respectively. The equality $\sigma_{|_{\Dho}} \fp = \fp \sigma$ actually shows that
\begin{equation}\label{hsigmafpmo}
H \sigma {\fp}\mo=H {\fp}\mo \sigma\,;
\end{equation}
so they represent the same hypermap.
\end{proof}

\begin{theorem}
If $\Wal (\cH) \cong \Wal (\cG)$, then $\cH \cong \cG$ or $\cH \cong \cD_{(0 1)} (\cG$).
\end{theorem}

\begin{proof}
If $\Wal (\cH) \cong \Wal (\cG)$ then $H\fw\mo=(G\fw\mo)^{g}$ for some $g\in\Delta$. \vskip4pt
\noi (i) $g\in\Dho$. Then $(G\fw\mo)^{g}=G^{{g\fw}}\fw\mo$, by Lemma \ref{generalprop}(1), and
then we have
\[H=H\fw\mo\fw=G^{{g\fw}}\fw\mo\fw=G^{{g\fw}};\]
that is, $\cH\isomorphic \cG$. \vskip4pt \noi (ii) $g\not\in\Dho$. Then $gR_0\in\Dho$ and
\[(G\fw\mo)^{g}=\left((G\fw\mo)^{gR_0}\right)^{R_0}=\left(G^{(gR_0)\fw}\fw\mo\right)^{R_0}=
G^{(gR_0)\fw}\sigma\fw\mo\,,\] using \eqref{hsigmafwmo}, where $\lambda=\iota^{R_0}$ and
$\sigma=(0\,1)\op$. Therefore
\[H=H\fw\mo\fw=G^{{(gR_0)\fw}}\sigma\fw\mo\fw=G^{{(gR_0)\fw}}\sigma,\]
which says that $\cH\isomorphic D_{\sigma}(\cG)$.
\end{proof}

\begin{theorem}
If $\Pin (\cH) \cong \Pin (\cG)$, then $\cH \cong \cG$.
\end{theorem}

\begin{proof}
As before, let $H$ and $G$ be hypermap-subgroups of $\cH$ and $\cG$. If $\Pin(\cH)\isomorphic
\Pin(\cG)$ then $H\fp\mo=(G\fp\mo)^{g}$ for some $g\in\Delta$. \vskip4pt \noi (i) If $g\in\Dho$
then, as before, $(G\fp\mo)^{g}=G^{g\fp}\fp\mo$ and then $H=G^{g\fp}$, showing that
$\cH\isomorphic \cG$. \vskip4pt \noi (ii) Suppose that $g\not\in\Dho$. As for $b\in\Dho$,
$(R_1R_2)^{R_0b}\fp=1\in H\nn G$ so that $(R_1R_2)^{R_0b}$ belongs to both $H\fp\mo$ and
$G\fp\mo$, for all $b\in\Dho$. Then (1) $R_1R_2\in (H\fp\mo)^{b\mo R_0}$ and (2) since
$(R_1R_2)^{R_0bg}\in (G\fp\mo)^g=H\fp\mo$, \, $R_1R_2\in (H\fp\mo)^{g\mo b\mo R_0}$. Since $b\mo
R_0$ runs all over $\Delta\m \Dho$ and $g\mo b\mo R_0$ runs all over $\Dho$, when $b\in\Dho$, then
$R_1R_2\in (H\fp\mo)^d$, for all $d\in\Delta$. This implies that all the hypervertices of
$\Pin(\cH)$ have valency 1. By a dual version of Lemma \ref{valencyone}, $\Pin(\cH)$ is a
``star''-like hypermap (see Figure \thefigcounter);

\ShowFigureVI{.6}{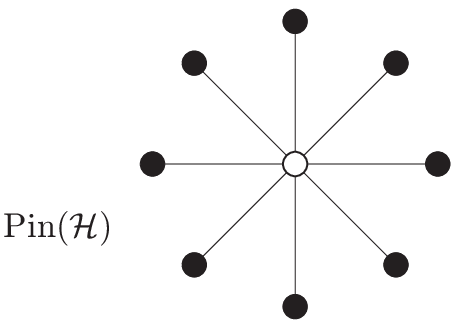}{$\Pin(\cH)=D_{(0\,2)} (\cD_n)$.} %
that is, $\Pin(\cH)=D_{(0\,2)} (\cD_n)$. Hence $\Pin(\cH)$ is a regular hypermap on the sphere with $n$ (even)
hypervertices. Thus $H\fp\mo$, as well as $(G\fp\mo)^g$, is normal in $\Delta$. Therefore, $H\fp\mo=G\fp\mo$ and hence
$H=G$.
\end{proof}

\noi The proof of the above theorem reveals the following information,

\begin{lemma} If $\Pin(\cH)$ is not isomorphic to $D_{(0\,2)} (\cD_{n})$ for any even $n$, then $\Pin(\cH)\isomorphic
\Pin(\cG)$ implies that $H\fp\mo=(G\fp\mo)^g$ for some $g\in\Dho$.
\end{lemma}

\subsection{Euler formula for bipartite-uniform hypermaps}

In this subsection we write the Euler characteristic of a bipartite-uniform hypermap in terms of
its bipartite-type. Let $\cH = (\Omega_{\cH}; h_0, h_1, h_2)$ be a bipartite-uniform hypermap with
Euler characteristic $\chi$, let $V$, $E$ and $F$ be the numbers of hypervertices, hyperedges and
hyperfaces of $\cH$, respectively, and let $V_1$ and $V_2 = V - V_1$ be the numbers of
hypervertices of the two $\Dho$-orbits in $\Omega_{\cH}$. By Lemma \ref{effh}, \,\, $\chi = V_1 +
V_2 + E + F - \frac{|\Omega_{\cH}|}{2}$. Let $(l_1, l_2; m; n)$ be the bipartite-type of $\cH$.
Then $V_1 = \frac{|\Omega_{\cH}|}{4 l_1}$, $V_2 = \frac{|\Omega_{\cH}|}{4 l_2}$, $E =
\frac{|\Omega_{\cH}|}{2m}$ and $F = \frac{|\Omega_{\cH}|}{2n}$. Replacing these values in the
above formula we get the following result:

\begin{lemma}[Euler formula for bipartite-uniform hypermaps]
If $\cH$ is a bipartite-uniform hypermap of bipartite-type $(l_1, l_2; m; n)$, then
\[\chi = \frac{|\Omega_{\cH}|}{2} \left( \frac{1}{2 l_1} + \frac{1}{2 l_2} + \frac{1}{m} + \frac{1}{n} - 1 \right).\]
\end{lemma}

\def\buH{\cK}
\def\hsbuH{K}
\def\closureDbuH{\nclosure{\buH}{\Delta}}
\def\coreDbuH{\core{\buH}{\Delta}}

\subsection{Spherical bipartite-uniform hypermaps}
In this subsection we classify the bipartite-uniform hypermaps $\buH$ on the sphere. The main
results were already given before; all we need now is to apply them directly to the sphere ($\chi
= 2$).

Let $\buH$ be a bipartite-uniform hypermap of bipartite-type $(l_1, l_2; m; n)$ on the sphere.
Then $\chi = 2 > 0$ and $\frac{1}{2 l_1} + \frac{1}{2 l_2} + \frac{1}{m} + \frac{1}{n} > 1$.
Suppose, without loss of generality, that $l_1 \leq l_2$ and $m \leq n$. Then

\begin{eqnarray*}
\frac{1}{l_1} + \frac{2}{m} \geq \frac{1}{2 l_1} + \frac{1}{2 l_2} + \frac{1}{m} + \frac{1}{n} > 1
&
\Rightarrow & \frac{1}{l_1} > \frac{1}{2} \quad \mbox{or} \quad \frac{2}{m} > \frac{1}{2} \\
& \sse & l_1 < 2 \quad \mbox{or} \quad m < 4 \\
& \sse & l_1 = 1 \quad \mbox{or} \quad m = 2 \\
& & \mbox{(since $m$ is even)}
\end{eqnarray*}

\noindent From this result and Theorems \ref{HWal} and \ref{HPin}, we deduce the following
theorem.

\begin{theorem}
If $\buH$ is a spherical bipartite-uniform hypermap, then $\buH \cong \Wal (\cR)$ or $\buH \cong
\Pin (\cR)$ for some spherical uniform hypermap $\cR$, unique up to isomorphism. Moreover, as
$\buH$ is bipartite-regular \ifif $\cR$ is regular, and on the sphere all uniform hypermaps are
regular, then all bipartite-uniform hypermaps on the sphere are bipartite-regular.
\end{theorem}

\begin{table}[ht]
{\scriptsize
\[\begin{array}{lllll|llll|llcl}
\mathbf{\#} & l_1 & l_2 & m & n & V_1 & V_2 & E & F & |\Omega| & \hskip8pt\buH \\
\hline
\mathbf{1} & 1 &  1 & 2n & 2n &  n &  n &  1 &  1 &  4n & \Pin (\cD_{(02)}(\cD_n)) \\
\mathbf{2} & 1 &  2 &  4 & 2n & 2n &  n &  n &  2 &  8n & \Pin (\cP_n) \\
\mathbf{3} & 1 &  2 &  6 &  6 & 12 &  6 &  4 &  4 &  48 & \Pin (\cD_{(01)} (\cT)) \\
\mathbf{4} & 1 &  2 &  6 &  8 & 24 & 12 &  8 &  6 &  96 & \Pin (\cD_{(01)} (\cC)) \\
\mathbf{5} & 1 &  2 &  6 & 10 & 60 & 30 & 20 & 12 & 240 & \Pin (\cD_{(01)} (\cD)) \\
\mathbf{6} & 1 &  3 &  4 &  6 & 12 &  4 &  6 &  4 &  48 & \Pin (\cT) \\
\mathbf{7} & 1 &  3 &  4 &  8 & 24 &  8 & 12 &  6 &  96 & \Pin (\cC) \\
\mathbf{8} & 1 &  3 &  4 & 10 & 60 & 20 & 30 & 12 & 240 & \Pin (\cD) \\
\mathbf{9} & 1 &  4 &  4 &  6 & 24 &  6 & 12 &  8 &  96 & \Pin (\cD_{(02)} (\cC)) \\
\mathbf{10} & 1 &  5 &  4 &  6 & 60 & 12 & 30 & 20 & 240 & \Pin (\cD_{(02)} (\cD)) \\
\mathbf{11} & 1 &  n &  2 & 2n &  n &  1 &  n &  1 &  4n & \Pin (\cD_{(12)} (\cD_n)) \\
\mathbf{12} & 1 &  n &  4 &  4 & 2n &  2 &  n &  n &  8n & \Pin (\cD_{(02)} (\cP_n)) \\
\mathbf{13} & 2 &  2 &  2 & 2n &  n &  n & 2n &  2 &  8n & \Wal (\cP_n) \\
\mathbf{14} & 2 &  3 &  2 &  6 &  6 &  4 & 12 &  4 &  48 & \Wal (\cT) \\
\mathbf{15} & 2 &  3 &  2 &  8 & 12 &  8 & 24 &  6 &  96 & \Wal (\cC) \\
\mathbf{16} & 2 &  3 &  2 & 10 & 30 & 20 & 60 & 12 & 240 & \Wal (\cD) \\
\mathbf{17} & 2 &  4 &  2 &  6 & 12 &  6 & 24 &  8 &  96 & \Wal (\cD_{(02)} (\cC)) \\
\mathbf{18} & 2 &  5 &  2 &  6 & 30 & 12 & 60 & 20 & 240 & \Wal (\cD_{(02)} (\cD)) \\
\mathbf{19} & 2 &  n &  2 &  4 &  n &  2 & 2n &  n &  8n & \Wal (\cD_{(02)} (\cP_n)) \\
\mathbf{20} & 3 &  3 &  2 &  4 &  4 &  4 & 12 &  6 &  48 & \Wal (\cD_{(12)} (\cT)) \\
\mathbf{21} & 3 &  4 &  2 &  4 &  8 &  6 & 24 & 12 &  96 & \Wal (\cD_{(12)} (\cC)) \\
\mathbf{22} & 3 &  5 &  2 &  4 & 20 & 12 & 60 & 30 & 240 & \Wal (\cD_{(12)} (\cD)) \\
\mathbf{23} & n &  n &  2 &  2 &  1 &  1 &  n &  n &  4n & \Wal (\cD_n) \\
  \hline
\end{array}\]
\caption{\label{biptable} The bipartite-regular hypermaps on the sphere.}}
\end{table}

Based on the knowledge of regular hypermaps on the sphere, we display in Table \ref{biptable} all
the possible values (up to duality) for the bipartite-type of the bipartite-regular hypermaps on
the sphere and the unique hypermap (up to isomorphism) with such a bipartite-type. Notice that the
map of bipartite-type $(1,n;2;2n)$ can be constructed from $\cD_n$ either via a Wal transformation
$\Wal (\cD_{(02)} (\cD_n))$ or via a Pin transformation $\Pin (\cD_{(12)} (\cD_n))$. Since $\Wal
(\cD_{(02)} (\cD_n))$ $\isomorphic$ $\Wal (\cD_{(12)} (\cD_n))$ these two constructions (Wal and
Pin) can actually be carried forward on the same hypermap $\cD_{(12)} (\cD_n)$. The Tetrahedron
$\cR=\cT$, which is self-dual, gives rise to $\Wal (\cD_{(0\;1)} (\cT)) = \Wal (\cT) = \Wal
(\cD_{(02)}(\cT))$.

\section{Irregularity and chirality}

We follow the same terminology and notations used in \cite{BCDirrh}. Let $\buH$ be a bipartite
(that is, $\Dho$-conservative) hypermap with hypermap-subgroup $\hsbuH<\Dho$. If $\buH$ is not
regular (that is, not $\De$-regular), then its closure cover $\closureDbuH$ is the largest regular
hypermap covered by $\buH$ and its covering core $\coreDbuH$ is the smallest regular hypermap
covering $\buH$. Hence we have two normal subgroups in $\De$, the normal closure $\hsbuH\sup{\De}$
containing $\hsbuH$, and the core $\hsbuH\sub{\De}$ contained in $\hsbuH$. Since
$\hsbuH\sub{\De}\normal \hsbuH$, although $\hsbuH$ may not be normal in $\hsbuH\sup{\De}$, we have
a group
\[\ligD(\buH)={\hsbuH/ \hsbuH\sub{\De}}\]
called the {\em lower-irregularity group} of $\buH$. Its size is the {\em lower-irregularity
index} and is denoted by $\liiD(\buH)$. The {\em upper-irregularity index}, denoted by
$\uiiD(\buH)$, is the index $|\hsbuH\sup{\De}:\hsbuH|$. If $\buH$ is bipartite-regular, then
$\hsbuH\normal \Dho$, and since $\hsbuH\sup{\De}$ is a subgroup of $\Dho$, $\hsbuH\normal
\hsbuH\sup{\De}$ and we have another group, the {\em upper-irregularity group}
\[\uigD(\buH)={\hsbuH\sup{\De}/ \hsbuH}.\]
Since the index of $\Dho$ in $\Delta$ is 2, the upper- and lower-irregularity groups are
isomorphic; so their upper- and lower-irregularity indices are equal ($\buH$ is irregularity
balanced). The common group $\uigD(\buH)\isomorphic \ligD(\buH)=\Upsilon$ is the {\em irregularity
group} of the bipartite-regular hypermap $\buH$ and the common value
$\uiiD(\buH)=\liiD(\buH)=\iota$ is its {\em irregularity index}. This has value 1 \ifif $\buH$ is
regular. Being bipartite-regular, $\buH$ is isomorphic to a regular $\Dho$-marked hypermap (see
\cite{Btrrh})
\[\cQ=(G,a,b,c,d)\isomorphic (\Dho/\hsbuH,\hsbuH A,\hsbuH B,\hsbuH C,\hsbuH D)\,,\]
where $\Dho=\gpg{A,B,C,D}\isomorphic C_2*C_2*C_2*C_2$ and $\hsbuH$ is the $\Dho$-hypermap subgroup
of $\cQ$ (and the hypermap subgroup of $\buH$). Here $G$ is the group generated by $a,b,c,d$. To
compute the irregularity group of $\buH$ we use:

\begin{lemma}\label{irregAB}
If $G$ has presentation $\langle a,b,c,d\st R=1\rangle$, where $R=\{R_1,\dots,R_k\}$ is a set of
relators $R_i=R_i(a,b,c,d)$ then $\uigD(\buH)=\< R^{R_0}\>\sup{G}$.
\end{lemma}
See \cite{BCDirrh} for the proof.

\def\nchsbuHD{\ncD{\hsbuH}}

\vskip5pt The definition of chirality given in \cite{Brch} is slightly different from that used in
\cite{BJNScgmh,BNchsg,BNchfh,BNjih}. If $\buH$ is bipartite ($\hsbuH<\Dho$), not necessarily
bipartite-regular, then $\buH$ is {\em $\Dho$-chiral}, or bipartite-chiral, if the normaliser
$N\sub{\De}(\hsbuH)$ of $\hsbuH$ in $\Delta$ is a subgroup of $\Dho$. In other words, $\buH$ is
$\Dho$-chiral if the group of automorphisms $Aut(\buH)\isomorphic N\sub{\De}(\hsbuH)/\hsbuH$
contains no ``symmetry'' besides $\Dho$.

Let $\buH$ be a $\Dho$-chiral hypermap. If $\buH$ is bipartite-regular ($\Dho$-regular), then
$\hsbuH\normal \Dho$ and so we have $N\sub{\De}(\hsbuH)=\Dho$. Thus $\buH$ is $\Dho$-chiral \ifif
$\hsbuH$ is not normal in $\De$; that is, \ifif $\buH$ is irregular. As $\Dho$ has index 2 in
$\Delta$, with transversal $\{1,R_0\}$, we have $\hsbuH=\hsbuH^{\<R_0\>}=\hsbuH
\hsbuH^{R_0}=\nchsbuHD$ \ifif $R_0\in N\sub{\De}(\hsbuH)$; that is, \ifif $\hsbuH R_0\in
Aut(\buH)$. Hence the upper-irregularity index $\uiiD$ gives a ``measure'' of ``how close'' $\buH$
is to having the ``symmetry'' $\hsbuH R_0$ outside $\Dho$. For this reason we also call the
upper-irregularity index (which coincides with the lower-irregularity index) the {\em
$\Dho$-chirality index} of the bipartite-regular $\buH$. This expresses how ``close'' $\buH$ is to
getting a ``symmetry'' outside $\Dho$, or in other words, how close it is for losing
$\Dho$-chirality.

The same happens to any normal subgroup $\Theta$ with index two in $\Delta$. In particular, for
$\Theta=\De^+$, the upper irregularity index (or simply the irregularity index) of a
$\De^+$-regular (that is, orientably regular) hypermap coincides with the $\De^+$-chirality index.
This explains the use of chirality index in place of irregularity index (of orientably regular
hypermaps) in the papers \cite{BJNScgmh,BNchsg,BNchfh,BNjih}. For more information and a general
definition of chirality group see \cite{Brch}.

\def\hsp{\hskip8pt}
\def\tasI{\!\!\!}
\def\tasF{\!\!\!\!}
\def\tasM{\!}
\def\tasMM{\!\!\!}
\def\tasMz{\!\!\!\!\!\!\!}
\def\tasMF{\!\!\!\!\!}
\def\tasMA{\!\!\!\!\!\!}
\def\ts{\!\!}
\def\tasIL{\!\!\!}
\def\tasFL{\!\!\!\!\!\!\!\!\!\!}

\def\HD#1{{\buH_{\De}}_{#1}}
\def\zs#1#2{\begin{array}{l} \tasI #1 \\ \tasI #2 \end{array}}

\begin{table}[ht]
{\scriptsize
\[\begin{array}{l|lll|llll|ll}
\tasIL\hsp\buH \tasF& \hsp\buH^{\Delta} \tasMz&\tasM l,m,n \tasMM&\tasM |\Omega| \tasMA&\tasM \hsp\buH_{\Delta} \tasMM&\tasM l,m,n \tasMM&\tasM |\Omega| \tasMM&\tasM \tx{genus} \tasMM&\tasM \,\iota \tasMF&\tasM \ig \tasFL\\
\hline %
\tasIL\Pin (\cD\sub{\ts(0 2)} (\cD\sub{n})) \tasF& \cD\sub{\ts(0 2)} (\cD\sub{2n}) \tasMz&\tasM 1,\!2n,\!2n
\tasMM&\tasM 4n \tasMA&\tasM \cD\sub{\ts(0 2)} (\cD\sub{2n}) \tasMM&\tasM 1,\!2n,\!2n \tasMM&\tasM 4n \tasMM&\tasM 0
\tasMM&\tasM 1 \tasMF&\tasM 1
\tasFL \\
\tasIL\Pin (\cP_n) \tasF& \tasI \left\{\zs{\cD\sub{\ts(0 2)} (\cD_4)}{\cD\sub{\ts(0 2)} (\cD_2)}\right. \tasMz&\tasM
\zs{1,4,4}{1,2,2} \tasMM&\tasM \zs{8}{4} \tasMA&\tasM \zs{\tasM \HD{2} }{\tasM \HD{3} } \tasMM&\tasM
\zs{2,4,2n}{2,4,2n} \tasMM&\tasM \zs{8n^2}{16n^2} \tasMM&\tasM \zs{\frac{(n-1)^2+1}{2}}{(n-1)^2} \tasMM&\tasM
\zs{n}{2n} \tasMF&\tasM \zs{D\sub{\frac{n}{2}}, \ \mbox{$n$ even}}{D_n, \ \mbox{$n$ odd}}
\tasFL \\
\tasIL\Pin (\cD\sub{\ts(0 1)} (\cT)) \tasF& \cD\sub{\ts(0 2)} (\cD_6) \tasMz&\tasM %
1,6,6 \tasMM&\tasM 12 \tasMM&\tasM \HD{4} \tasMM&\tasM 2,6,6 \tasMM&\tasM 192 \tasMM&\tasM 9 \tasMM&\tasM 4
\tasMF&\tasM V_4
\tasFL \\
\tasIL\Pin (\cD\sub{\ts(0 1)} (\cC)) \tasF& \cD\sub{\ts(0 2)} (\cD_2) \tasMz&\tasM %
1,2,2 \tasMM&\tasM 4 \tasMA&\tasM \HD{5} \tasMM&\tasM 2,6,8 \tasMM&\tasM 2304 \tasMM&\tasM 121 \tasMM&\tasM 24
\tasMF&\tasM S_4
\tasFL \\
\tasIL\Pin (\cD\sub{\ts(0 1)} (\cD)) \tasF& \cD\sub{\ts(0 2)} (\cD_2) \tasMz&\tasM %
1,2,2 \tasMM&\tasM 4 \tasMA&\tasM \HD{6} \tasMM&\tasM 2,6,10 \tasMM&\tasM 14400 \tasMM&\tasM 841 \tasMM&\tasM 60
\tasMF&\tasM A_5
\tasFL \\
\tasIL\Pin (\cT) \tasF& \cD\sub{\ts(0 2)} (\cD_2) \tasMz&\tasM %
1,2,2 \tasMM&\tasM 4 \tasMA&\tasM \HD{7} \tasMM&\tasM 3,4,6 \tasMM&\tasM 576 \tasMM&\tasM 37 \tasMM&\tasM 12
\tasMF&\tasM A_4
\tasFL \\
\tasIL\Pin (\cC) \tasF& \cD\sub{\ts(0 2)} (\cD_4) \tasMz&\tasM %
1,4,4 \tasMM&\tasM 8 \tasMA&\tasM \HD{8} \tasMM&\tasM 3,4,8 \tasMM&\tasM 1152 \tasMM&\tasM 85 \tasMM&\tasM 12
\tasMF&\tasM A_4
\tasFL \\
\tasIL\Pin (\cD) \tasF& \cD\sub{\ts(0 2)} (\cD_2) \tasMz&\tasM %
1,2,2 \tasMM&\tasM 4 \tasMA&\tasM \HD{9} \tasMM&\tasM 3,4,10 \tasMM&\tasM 14400 \tasMM&\tasM 1141 \tasMM&\tasM 60
\tasMF&\tasM A_5
\tasFL \\
\tasIL\Pin (\cD\sub{\ts(0 2)} (\cC)) \tasF& \cD\sub{\ts(0 2)} (\cD_2) \tasMz&\tasM %
1,2,2 \tasMM&\tasM 4 \tasMA&\tasM \HD{10} \tasMM&\tasM 4,4,6 \tasMM&\tasM 2304 \tasMM&\tasM 193 \tasMM&\tasM 24
\tasMF&\tasM S_4
\tasFL \\
\tasIL\Pin (\cD\sub{\ts(0 2)} (\cD)) \tasF& \cD\sub{\ts(0 2)} (\cD_2) \tasMz&\tasM %
1,2,2 \tasMM&\tasM 4 \tasMA&\tasM \HD{11} \tasMM&\tasM 5,4,6 \tasMM&\tasM 14400 \tasMM&\tasM 1381 \tasMM&\tasM 60
\tasMF&\tasM A_5
\tasFL \\
\tasIL\Pin (\cD\sub{\ts(1 2)} (\cD_n)) \tasF& \cD\sub{\ts(0 2)} (\cD_2) \tasMz&\tasM %
1,2,2 \tasMM&\tasM 4 \tasMA&\tasM \HD{12} \tasMM&\tasM n,2,2n \tasMM&\tasM 4n^2 \tasMM&\tasM \frac{(n-1)(n-2)}{2}
\tasMM&\tasM n \tasMF&\tasM C_n
\tasFL \\
\tasIL\Pin (\cD\sub{\ts(0 2)} (\cP_n)) \tasF& \cD\sub{\ts(0 2)} (\cD_4) \tasMz&\tasM %
1,4,4 \tasMM&\tasM 8 \tasMA&\tasM \HD{13} \tasMM&\tasM n,4,4 \tasMM&\tasM 8n^2 \tasMM&\tasM (n-1)^2 \tasMM&\tasM n
\tasMF&\tasM C_n
\tasFL \\
\tasIL\Wal (\cP_n) \tasF& \cP\sub{2n} \tasMz&\tasM %
2,2,2n \tasMM&\tasM 8n \tasMA&\tasM \cP\sub{2n} \tasMM&\tasM 2,2,2n \tasMM&\tasM 8n \tasMM&\tasM 0 \tasMM&\tasM 1
\tasMF&\tasM 1
\tasFL \\
\tasIL\Wal (\cT) \tasF& \cD\sub{\ts(0 2)} (\cD_2) \tasMz&\tasM %
1,2,2 \tasMM&\tasM 4 \tasMA&\tasM \HD{15} \tasMM&\tasM 6,2,6 \tasMM&\tasM 576 \tasMM&\tasM 25 \tasMM&\tasM 12
\tasMF&\tasM A_4
\tasFL \\
\tasIL\Wal (\cC) \tasF& \cD\sub{\ts(0 2)} (\cD_2) \tasMz&\tasM %
1,2,2 \tasMM&\tasM 4 \tasMA&\tasM \HD{16} \tasMM&\tasM 6,2,8 \tasMM&\tasM 2304 \tasMM&\tasM 121 \tasMM&\tasM 24
\tasMF&\tasM S_4
\tasFL \\
\tasIL\Wal (\cD) \tasF& \cD\sub{\ts(0 2)} (\cD_2) \tasMz&\tasM %
1,2,2 \tasMM&\tasM 4 \tasMA&\tasM \HD{17} \tasMM&\tasM 6,2,10 \tasMM&\tasM 14400 \tasMM&\tasM 841 \tasMM&\tasM 60
\tasMF&\tasM A_5
\tasFL \\
\tasIL\Wal (\cD\sub{\ts(0 2)} (\cC)) \tasF& \cP_6 \tasMz&\tasM %
2,2,6 \tasMM&\tasM 24 \tasMA&\tasM \HD{18} \tasMM&\tasM 4,2,6 \tasMM&\tasM 384 \tasMM&\tasM 9 \tasMM&\tasM 4
\tasMF&\tasM V_4
\tasFL \\
\tasIL\Wal (\cD\sub{\ts(0 2)} (\cD)) \tasF& \cD\sub{\ts(0 2)} (\cD_2) \tasMz&\tasM %
1,2,2 \tasMM&\tasM 4 \tasMA&\tasM \HD{19} \tasMM&\tasM 10,2,6 \tasMM&\tasM 14400 \tasMM&\tasM 841 \tasMM&\tasM 60
\tasMF&\tasM A_5
\tasFL \\
\tasIL\Wal (\cD\sub{\ts(0 2)} (\cP_n)) \tasF& \tasI \left\{\zs{\cP_4}{\cD\sub{\ts(0 2)} (\cD_2)}\right. \tasMz&\tasM
\zs{2,2,4}{1,2,2} \tasMM&\tasM \zs{16}{4} \tasMA&\tasM \zs{\tasM \HD{20} }{\tasM \HD{21} } \tasMM&\tasM
\zs{n,2,4}{2n,2,4} \tasMM&\tasM \zs{4n^2}{16n^2} \tasMM&\tasM \zs{\frac{(n-2)^2}{4}}{(n-1)^2} \tasMM&\tasM
\zs{\frac{n}{2}}{2n} \tasMF&\tasM \zs{C\sub{\frac{n}{2}}, \ \mbox{$n$ even}}{D_n, \ \mbox{$n$ odd}}
\tasFL \\
\tasIL\Wal (\cD\sub{\ts(1 2)} (\cT)) \tasF& \cC \tasMz&\tasM %
3,2,4 \tasMM&\tasM 48 \tasMA&\tasM \cC \tasMM&\tasM 3,2,4 \tasMM&\tasM 48 \tasMM&\tasM 0 \tasMM&\tasM 1 \tasMF&\tasM 1
\tasFL \\
\tasIL\Wal (\cD\sub{\ts(1 2)} (\cC)) \tasF& \cD\sub{\ts(0 2)} (\cD_2) \tasMz&\tasM %
1,2,2 \tasMM&\tasM 4 \tasMA&\tasM \HD{23} \tasMM&\tasM 12,2,4 \tasMM&\tasM 2304 \tasMM&\tasM 97 \tasMM&\tasM 24
\tasMF&\tasM S_4
\tasFL \\
\tasIL\Wal (\cD\sub{\ts(1 2)} (\cD)) \tasF& \cD\sub{\ts(0 2)} (\cD_2) \tasMz&\tasM %
1,2,2 \tasMM&\tasM 4 \tasMA&\tasM \HD{24} \tasMM&\tasM 15,2,4 \tasMM&\tasM 14400 \tasMM&\tasM 661 \tasMM&\tasM 60
\tasMF&\tasM A_5
\tasFL \\
\tasIL\Wal (\cD_n) \tasF& \cD\sub{\ts(0 2)} (\cP_n) \tasMz&\tasM %
n,2,2 \tasMM&\tasM 4n \tasMA&\tasM \cD\sub{\ts(0 2)} (\cP_n) \tasMM&\tasM n,2,2 \tasMM&\tasM 4n \tasMM&\tasM 0 \tasMM&\tasM 1 \tasMF&\tasM 1\tasFL \\
\hline %
\end{array}\]
\caption{\label{xtable} $\buH$, $\buH^{\Delta}$ and $\buH_{\Delta}$.}}
\end{table}

\vskip5pt \noi {\bf Computing the irregularity group $\mathbf{\ig}$.} \vskip5pt

\noi Let $\buH=\tfH=\Wal(\cH)$ or $\Pin(\cH)$ conform $\vphi=\fw$ or $\fp$, respectively, and let
$H$ be the hypermap subgroup of a regular hypermap $\cH$ of type $(l;m;n)$. The inverse image
$\hsbuH=H\vphi\mo$ is the hypermap subgroup of $\buH$. The lower-irregularity index of $\buH$,
$\ligD(\buH) = {\hsbuH/ \hsbuH\sub{\De}}$, is isomorphic to its upper-irregularity group
$\uigD(\buH) = {{\hsbuH}\sup{\De}/ \hsbuH}$, a subgroup of the $\Dho$-monodromy group
$G=\Dho/\hsbuH$ of $\buH$. This common group, the irregularity group $\ig$, can be computed in the
following way. According to Theorem \ref{Dhomonodromy}, the group $G\isomorphic \Mon(\cH)$ is a
known group (see Table \ref{sphereunif}). Being $G$ the $\Dho$-monodromy group of a
bipartite-regular hypermap, using $\varphi$ we can rewrite $G$ in the following form
\[G=\langle a,b,c,d\st a^2=b^2=c^2=d^2=1, R=1\rangle\,,\]
such that $a^{R_0}=c$, $b^{R_0}=d$, $c^{R_0}=a$ and $d^{R_0}=b$; \, $R$ stands for a set of
relators on $a,b,c,d$. By Lemma \ref{irregAB},
\[\ig=\langle R^{R_0}\rangle\sup{G}\]
is the closure subgroup of $R^{R_0}$ in $G$. This calculation is easily performed and the results
for $\ig(\buH)$ can be seen under the last column of Table \ref{xtable}. For an example of how
this calculation is carried out see Theorem \ref{bpchiralhpsoneorsurface}.

\vskip5pt

However, since $\vphi=\fw$ or $\fp$ sends generators of $\Dho$ of odd length in $\De$ to
generators of $\De$ of odd length in $\De$, we have necessarily $\Dp\vphi\mo=\Dpho$, where
$\Dpho=\Dp \nn \Dho$\oop{. In fact, $\Dp\normal_2\De$ implies $\Dp\vphi\mo\normal_2\Dho$. Now
$\Dpho\normal_2\Dho$. Since $\Dpho\vphi\cc \Dp$ then $\Dpho\cc \Dpho\vphi \vphi\mo\cc \Dp\vphi\mo$
and so $\Dp\vphi\mo=\Dpho$}. Since $H\normal \Dp$, then $K\sup{\De}/K\normal \Dpho /
K=\Dp\vphi\mo/H\vphi\mo \isomorphic \Dp / H=Aut^{+}(\cH)$\oop{, where $Aut^{+}(\cH)$ is the group
of the orientation preserving automorphisms of $\cH$}; that is $\ig=K\sup{\De}/K$ is a normal
subgroup of $Aut^{+}(\cH)$.

\vskip4pt

Let $A=R_1$, $B=R_2$, $C=R_1^{R_0}$ and $D=R_2^{R_0}$. Then $\Dho=\gpg{A,B,C,D}$.

\vskip4pt \noindent (1) If $\buH=Wal(\cH)$, then $K=\gpg{BD, (AB)^l, (DC)^m, (CA)^n}^{\Dho}$, so
$K^{\De}=\gpg{BD, (AB)^l, (DC)^m, (CA)^n}^{\De}$. Let $d=\gcd(l,m)$. Since
$(AB)^m=((DC)^{-m})^{R_0}$ and $(DC)^l=((AB)^{-l})^{R_0}$ then $(AB)^d$ and $(DC)^d$ also belong
to $K^{\De}$. Hence if $d=1$, then $AB$ and $DC$ belong to $K^{\De}$ and so $K^{\De}=\Dpho$.
Therefore $\ig=Aut^{+}(\cH)$ when $d=1$.

\vskip4pt \noindent (2) If $\buH=Pin(\cH)$, then $K=\gpg{CD, (AB)^l, (BD)^m, (CA)^n}^{\Dho}$, so
$K^{\De}=\gpg{CD, (AB)^l, (BD)^m, (CA)^n}^{\De}$. Let $d=\gcd(m,n)$. Since $K^{\De}D=K^{\De}C$,
$K^{\De}(CA)^m=K^{\De}(DA)^m=(K^{\De}(BC)^m)^{R_0}=(K^{\De}(BD)^m)^{R_0}=K^{\De}$ and so
$(CA)^m\in K^{\De}$. Similarly, $(BD)^n\in K^{\De}$. Hence if $d=1$, then $K^{\De}=\Dpho$ and
consequently $\ig=Aut^{+}(\cH)$. \vskip5pt

Therefore the general calculations mentioned above only need to be carried out for the cases where
$d\neq 1$, namely cases 2 (with $n$ even), 3, 7, 12, 17 and 19 (with $n$ even).

\vskip7pt \noi {\bf Computing the closure cover $\mathbf{\ccover{\buH}}$.} \vskip5pt \noi Once the
irregularity index is calculated, it is an easy task to compute the closure cover $\ccover{\buH}$
of $\buH=\tfH$\!\!\!, simply because the genus of the closure cover is zero and in the sphere the
type determines uniquely a uniform (or regular) hypermap. Let $(l;m;n)$ be the type of the closure
cover $\ccover{\buH}$ and let $(r,s;u;v)$ be the bipartite-type of the spherical bipartite-regular
hypermap $\buH$. The number of flags $|\Omega_{\ccover{\buH}}|$ of the closure cover must divide
the number of flags $|\Omega_{\buH}|$ of $\buH$. Also $l$ divides $gcd(r,s)$, $m$ divides $u$ and
$n$ divides $v$. The greatest possible values for $l$, $m$ and $n$ are $gcd(r,s)$, $u$ and $v$,
respectively. Moreover, when $gcd(r,s)=1$ we must have $l=1$ in which case $m=n$ and the greatest
possible values are achieved for $m=n=gcd(u,v)$. Since $\ccover{\buH}$ is a regular hypermap on
the sphere and is determined by $l$, $m$ and $n$, we must check if in each case the above choice
of $l$, $m$ and $n$ give rise to a spherical type (cf. Table \ref{sphereunif}). If not we choose
the second greatest, the third greatest and so forth. For each bipartite-regular hypermap $\buH$
in Table \ref{biptable}, where $(l_1,l_2;m;n)$ is our $(r,s;u;v)$, taking the greatest values for
the triple $(l,m,n)$ we get a spherical type. To check if such triple determines a hypermap
covered by $\buH$ we take a half-turn in the middle of each hyperedge of $\buH$; these half-turns
determine a covering \, $\buH\mapsto \ccover{\buH}$. The results can be seen in Table
\ref{xtable}.

\vskip7pt
\noi {\bf Computing the covering core $\ccore{\buH}$.} \vskip5pt \noi The covering core
is already computed since we know its monodromy group
\[Mon(\ccore{\buH})=Mon(\buH)\]
and their canonical generators. Feeding these parameters in GAP \cite{GAP4}, for example, we get
the rest of the information shown in the Table \ref{xtable}. In this table we observe two isolated
maps (not in families) with less then 100 edges, the map $D\sub{(1\,2)}(\HD{4})$ with 48 edges and
Petrie path of length 4, and $\HD{18}$ with 96 edges and Petrie path of length 6. In
\cite{Wthesis}, where we can find a good list of regular maps up to 100 edges (although the list
is guaranteed to be complete only up to 49), these maps are $P(70)$ and $DP(190)$ on pages 144 and
181 respectively. These can be consulted in the recently created Census of orientably-regular maps
\cite{Wcrm}.

\section{Final comments}

\noi Looking back at Table \ref{xtable} we clearly read the following extra result,

\begin{theorem}
The irregularity (or chirality) index of a bipartite-regular hypermap can be any positive integer
number. Moreover, cyclic groups and dihedral groups are irregularity groups of bipartite-regular
hypermaps.
\end{theorem}

\noi Using the $Pin$ and $Wal$ transformations we can say a little more,

\begin{theorem}\label{bpchiralhpsoneorsurface}
On each orientable surface of genus $g$ there are $\Dho$-chiral hypermaps (that is, irregular
bipartite-regular hypermaps) with irregularity indices $2g+1$, $4g+2$ and $4g$.
\end{theorem}
\begin{proof}
Just take the $Pin(\cM_k)$ and the $Wal(\cM_k)$ constructions over the one-face regular map
$\cM_k$ formed from a single $2k$-gon by identifying opposite edges orientably. The map $\cM_k$
has type $(k;2;2k)$ or $(2k;2;2k)$ according as $k$ is odd or even. The monodromy group of $\cM_k$
is the dihedral group $D_{2k}$ generated by the involutions $r_0$, $r_1$ and $r_2$ subject to the
relations $(r_0r_1)^{2k}=1$ and $r_2=r_0(r_1r_0)^k$. The genus of $\cM_k$ is $k-1\over 2$ if $k$
is odd and $k\over 2$ otherwise. Hence each orientable surface of genus $g$ supports two maps
$\cM_k$, one for $k$ odd and another for $k$ even. Note that $\cM_k$ has 1 or 2 vertices according
as $k$ is even or odd.

\ShowFigureVE{.3}{MkPinMkWalMk}{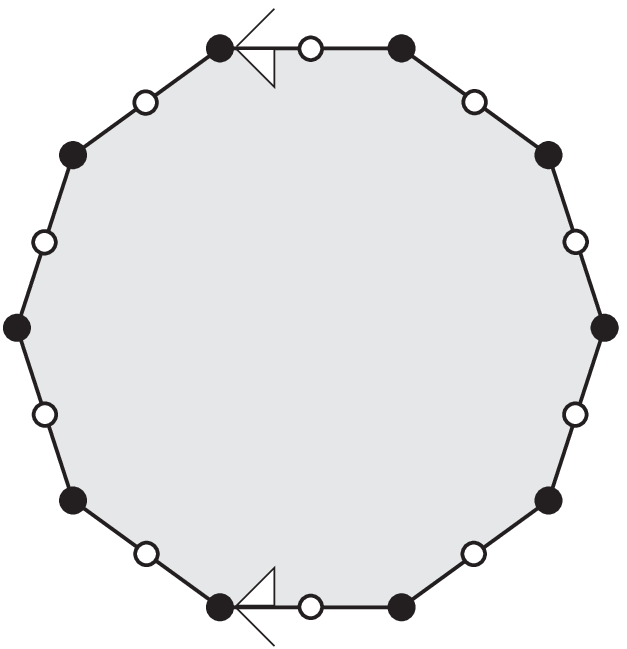}{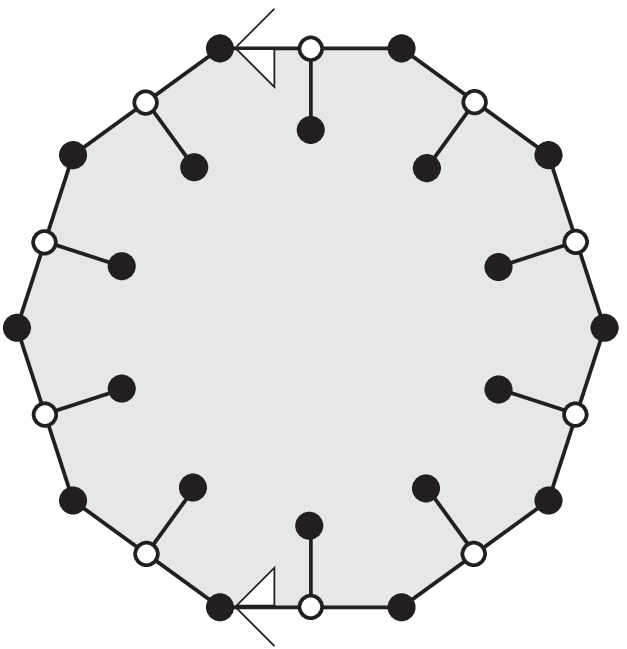}{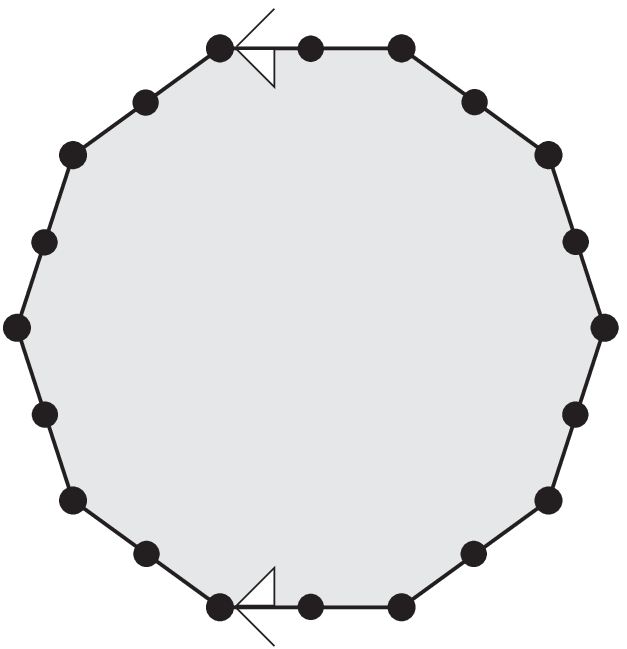}{(a) The $\cM_k$ map (opposite edges identified orientably).\\
(b) \, $Pin(\cM_k)$. \,\, (c)\, $Wal(\cM_k)$.}

\noindent The bipartite-regular hypermap $Pin(\cM_k)$ has bipartite-type $(1,k;4;4k)$ if $k$ is
odd and $(1,2k;4;4k)$ otherwise. The bipartite-regular map $Wal(\cM_k)$ has type $(2,k;2;4k)$ or
$(2,2k;2;4k)$ according as $k$ is odd or even. Let $H$ be the hypermap subgroup of $\cM_k$ and
$\hsbuH=H\varphi\mo$, where $\varphi=\fp$ or $\fw$.\oop{Recall that $\Dho=\gpg{A,B,C,D}\isomorphic
C_2*C_2*C_2*C_2$.}

\vskip4pt

\noindent (1) {\sl The hypermap $Pin(\cM_k)$}. The epimorphism $\fp$ induces an isomorphism
$G=\Dho/\hsbuH\para \Delta/H$, mapping $a\mapsto r_1$, $b\mapsto r_2$, $c\mapsto r_0$ and
$d\mapsto r_0$. That is, $c=d$ in the $\Dho$-monodromy group of $Pin(\cM_k)$. With the help of
$\fp$ we rewrite $Mon(\cM_k)$ in function of $a$, $b$, $c$ and $d$ to get the $\Dho$-monodromy
group
\[G=\gpg{a,b,c,d\st a^2=b^2=c^2=d^2=1, c=d, (ca)^{2k}=1, b=c(ac)^k}\,.\]
In this case $R=\{cd^{-1}, (ac)^{2k}, c(ac)^kb^{-1}\}$ and the irregularity group of $Pin(\cM_k)$
is the normal closure of $R^{R_0}$ in $G$; thus
\[\ig=\nclosureG{ab^{-1},(ca)^{2k},a(ca)^kd^{-1}}=\nclosureG{ab}=\gpg{ab}.\]
Since $ab=(ac)^{k+1}$, this group has size $k$ if $k$ is odd and size $2k$ otherwise. Hence
$Pin(\cM_k)$ has irregularity index $\ii=k=2g+1$, for $k$ odd, and $\ii=2k=4g$ for $k$ even.

\vskip4pt

\noindent (2) {\em The map $Wal(\cM_k)$}. Proceeding similarly we obtain
\[\begin{array}{rcl}
    G & = & \Dho-Mon(Wal(\cM_k)) \\
      & = & \gpg{a,b,c,d\st a^2=b^2=c^2=d^2=1, b=d, (ca)^{2k}=1, b=c(ac)^k}
  \end{array}\]
and irregularity group $\ig=\gpg{ac}=C_{2k}$ cyclic, giving rise to irregularity indices
$\ii=2k=4g+2$ when $k$ is odd and $\ii=2k=4g$ when $k$ is even.
\end{proof}

For non-orientable surfaces we cannot answer affirmatively since to obtain $\Dho$-chiral hypermaps
the $Pin$ and $Wal$ constructions need regular hypermaps and we know that there are none on the
non-orientable surfaces with negative characteristic 0, 1, 16, 22, 25, 37, and 46
\cite{BJrhgz,WBswnrh}.

\def\bookS#1#2#3#4#5{\bibitem{#1} #2, ``{\rm #3}'', #4, #5.}
\def\bookG#1#2#3#4#5#6{\bibitem{#1} #2, ``{\rm #3}'', #4, #5, #6.}  
\def\bookF#1#2#3#4#5#6#7{\bibitem{#1} #2, ``{\rm #3}'', {\sl #4}, #5, #6, #7.}  

\def\LectNotesF#1#2#3#4#5#6#7{\bibitem{#1} {\rm #2}, {\sl #3}, {\rm #4}, #5 (#6)(#7).} 
\def\LectNotesB#1#2#3#4#5#6#7#8{\bibitem{#1} {\rm #2}, {\sl #3}, {\rm #4}, #5, #6, #7, #8.} 

\def\paperS#1#2#3#4#5{\bibitem{#1} {\rm #2}, {\sl #3}, {\rm #4}, #5.}
\def\paperG#1#2#3#4#5#6{\bibitem{#1} {\rm #2}, {\sl #3}, {\rm #4}, #5, #6.}   
\def\paperF#1#2#3#4#5#6#7{\bibitem{#1} {\rm #2}, {\sl #3}, {\rm #4}, {\bf #5} (#6) #7.}   
\def\paperB#1#2#3#4#5#6#7#8{\bibitem{#1} {\rm #2}, {\sl #3}, {\rm #4}, (#5) {\bf #6} (#7) #8.} 
\def\papersub#1#2#3{\bibitem{#1} {\rm #2}, {\sl #3}, submitted.}
\def\paperappearing#1#2#3#4{\bibitem{#1} {\rm #2}, {\sl #3}, due to appear in the #4.}

\def\preprintS#1#2#3#4#5{\bibitem{#1} {\rm #2}, {\sl #3}, {#4}, #5.}
\def\proceedingsF#1#2#3#4#5#6#7{\bibitem{#1} {\rm #2}, {\sl #3}, {\rm #4}, {#5}, #6, #7.}   
\def\census#1#2#3#4{\bibitem{#1} {\rm #2}, {\rm #3}, {\rm #4}.}
\def\thesisS#1#2#3#4#5{\bibitem{#1} #2, {\sl #3}, ``{\rm Ph.D. Thesis}'', #4, #5.}
\def\thesisG#1#2#3#4#5#6{\bibitem{#1} #2, {\sl #3}, ``{\rm #4}'', #5, #6.}

\end{document}